\documentclass[twoside,leqno,twocolumn]{article}

\usepackage[letterpaper]{geometry}

\usepackage{ltexpprt}
\usepackage{hyperref}
\usepackage{amsfonts,amsmath}  
\usepackage{bm}
\usepackage{algorithm, algorithmic}
\usepackage{subfig,graphicx}

\newcommand{\KL}{K{\L}~}

\begin{document}

\newcommand\relatedversion{}

\title{\Large 
A Momentum Accelerated Adaptive Cubic Regularization Method for Nonconvex 
Optimization}
\author{Yihang Gao\thanks{Department of Mathematics, The University of Hong Kong, Pokfulam, Hong Kong SAR.}
\and Michael K. Ng\footnotemark[1]}

\date{}

\maketitle







\begin{abstract} \small\baselineskip=9pt 
The cubic regularization method (CR) 
and its adaptive version (ARC)
are popular Newton-type methods in solving unconstrained non-convex optimization problems, due to its global convergence to local minima under mild conditions. 
The main aim of this paper is to  
develop a momentum accelerated adaptive cubic regularization method (ARCm)
to improve the convergent performance.
With the proper choice of momentum step size, we show the  
global convergence of ARCm and the local convergence can also be guaranteed under the \KL property. Such global and local convergence can also be established when 
inexact solvers with low computational costs are employed in the iteration procedure. 
Numerical results for non-convex logistic regression and robust linear regression models 
are reported to demonstrate that
the proposed ARCm significantly outperforms state-of-the-art cubic regularization methods (e.g., CR, momentum-based CR, ARC) and the trust region method. In particular, 
the number of iterations required by ARCm is less than 
10\% to 50\% required by the most competitive method (ARC) in the experiments.
\end{abstract}


\section{Introduction}
Most machine learning tasks involve solving challenging (non-convex) optimization problems
\begin{equation}
\label{eq_intro_1}
    \min_{\bm{x} \in \mathbb{R}^{d}} f(\bm{x}).
\end{equation}
Second-order critical points are usually the preferred solutions for (\ref{eq_intro_1}) since many machine learning problems have been shown to have only strict saddle points and global minima without spurious local minima \cite{ge2016matrix,sun2016geometrical}. Second-order methods that exploit Hessian information have been proposed to escape saddle points \cite{nesterov2006cubic,conn2000trust} and enjoy faster local convergence \cite{li2004regularized,yue2019quadratic,zhou2018convergence}. Cubic regularization method (CR), first proposed by Griewank \cite{griewank1981modification}, and later independently by Nesterov and Polyak \cite{nesterov2006cubic}, and Weiser et al. \cite{weiser2007affine}, one of the second-order methods, is widely applied in solving inverse problems \cite{cartis2011adaptive,jiang2021accelerated}, regression models \cite{yue2019quadratic,zhou2018stochastic} as well as minimax problems \cite{huang2022cubic}. 

The vanilla CR method is formulated as 
\begin{equation*}
\begin{split}
    & \bm{s}_k \in \arg \min_{\bm{s}} \nabla f(\bm{x}_k)^{\top} \bm{s} + \frac{1}{2} \bm{s}^{\top} \nabla^{2} f(\bm{x}_k) \bm{s} + \frac{M_{k}}{6} \|\bm{s}\|_2^3,\\
    & \bm{x}_{k+1} =  \bm{x}_{k} + \bm{s}_k,
\end{split}
\end{equation*}
where $M_{k} := M$ is a fixed pre-defined constant for all iterations. In recent decades, various techniques are adopted to improve the performance of CR. Wang et al. \cite{wang2020cubic} accelerated the vanilla CR by a momentum term (CRm), 
where the method works well mainly by enlarging the step size of $\bm{s}_{k}$ when the cubic penalty parameter $M$ is over-estimated.

Cartis et al. \cite{cartis2011adaptive} proposed adaptive cubic regularization method (ARC), which adaptively assigns $M_{k}$ based on the quality of the step $\bm{s}_{k}$, as an analogy to the trust region method (TR) \cite{conn2000trust}. Stochastic (subsampled) ARC (e.g., Kohler et al. \cite{kohler2017sub} and Zhou et al. \cite{zhou2018stochastic} etc.) were developed to reduce the overall computation, where the gradient $\nabla f(\bm{x}_{k})$ and the Hessian $\nabla^2 f(\bm{x}_{k})$ are inexactly evaluated.  To the best of our knowledge, ARC behaves better in most of the applications compared with CR methods that fix $M$ (e.g., vanilla CR and CRm). A natural question is how to further accelerate ARC with little extra cost in each iteration. 

\begin{algorithm*}[t!]
 \caption{Adaptive cubic regularization with momentum (ARCm)}
 \label{alg_arcm}
 \begin{algorithmic}[1]
 \renewcommand{\algorithmicrequire}{\textbf{Input: $\bm{x}_{0}$, $\bm{v}_{-1} = \bm{0}$, $ \gamma_{1}>1 \geq \gamma_{2} > \gamma_{3}>0$, $1>\eta_{2} > \eta_{1}>0$, $\sigma_{0} > \sigma_{\text{min}}>0$, $1 > \tau, \beta >0$ and $\alpha_1, \alpha_2 > 0$}}
 \renewcommand{\algorithmicensure}{\textbf{Output: $\{\bm{x}_{k}\}_{k=0}^{T}$}}
 \REQUIRE in
 \ENSURE  out
  \FOR {$k = 0$ to $T-1$}
  \STATE Solve the cubic subproblem:
  \begin{equation*}
      \bm{s}_k \in \arg \min_{\bm{s}} m_k(\bm{s}) := \arg \min_{\bm{s}} f(\bm{x}_k) + \nabla f(\bm{x}_k)^{\top} \bm{s} + \frac{1}{2} \bm{s}^{\top} \nabla^{2} f(\bm{x}_k) \bm{s} + \frac{\sigma_{k}}{6} \|\bm{s}\|_2^3.
  \end{equation*}
  
  \STATE Compute $f(\bm{x}_k + \bm{s}_k)$ and 
      $\rho_{k} = {\displaystyle 
      \frac{f(\bm{x}_k) - f(\bm{x}_k + \bm{s}_k)}{f(\bm{x}_k) - m_k(\bm{s}_{k})}}$.
  
  \IF{$\rho_k > \eta_1$ (successful update)}
  \STATE $\bm{y}_{k+1} = \bm{x}_k + \bm{s}_k$
  \STATE \textbf{Momentum step:}\\ Select $\beta_{k} \in \left[0, \min \left(\tau, \alpha_1 \|\bm{s}_k\|_2, \alpha_2 \|\bm{s}_k\|_2^2 \right)\right]$ such that $f(\bm{z}_{k+1}) \leq f(\bm{y}_{k+1})$ with 
  $\bm{v}_{k} = \beta_{k} \cdot \bm{v}_{k-1}+\bm{s}_{k}~\text{and}~\bm{z}_{k+1} = \bm{x}_{k} + \bm{v}_{k}$.
  \STATE $\bm{x}_{k+1} = \bm{z}_{k+1}$ 
  \IF{$\rho_{k} > \eta_2$ (very successful update)}
  \STATE $\sigma_{k+1} = \max(\sigma_{\text{min}}, \gamma_3 \cdot \sigma_{k})$
  \ELSE
  \STATE $\sigma_{k+1} = \gamma_2 \cdot \sigma_{k}$
  \ENDIF
  \ELSE 
  \STATE (unsuccessful update)\\ $\bm{x}_{k+1}= \bm{x}_{k}$ and $\bm{v}_{k} = \bm{v}_{k-1}$
  \STATE $\sigma_{k+1} = \gamma_1 \cdot \sigma_{k}$
  \ENDIF
  \ENDFOR
 \end{algorithmic}
 \end{algorithm*}

In this paper, we investigate a momentum-accelerated adaptive cubic regularization method (ARCm). Here are our contributions and the outline of the paper.
\begin{itemize}
    \item We develop a momentum-accelerated adaptive cubic regularization method (ARCm). We adopt a general scheme for momentum, which is more suitable for ARC than CR (see in Theorem \ref{theorem_0}). The extra computation is cheap and ignorant as both the momentum and its step size are solely based on $\{\bm{s}_{k}\}$ and are free of gradient or Hessian evaluations. 

    \item With the proper setting of step size for momentum (see in the Algorithm \ref{alg_arcm}), the global convergence of ARCm to second-order critical points is satisfied (see in Theorem \ref{theorem_1}). We further show that the proposed ARCm enjoys local convergence under the \KL property (see in Theorem \ref{theorem_2}), which is one of the advantages of second-order methods over first-order methods.

    \item We also study the global convergence of ARCm with the inexact cubic regularized subproblems (CRS) solutions as we usually approximately solve CRS in practice (see in Theorem \ref{theorem_3}). The local convergence is preserved if the error of CRS is decreasing with $\|\bm{s}_{k}\|_2^3$ (see in Corollary \ref{corollary_1}) but may be destructed otherwise. 

    \item We conduct experiments in solving high-dimensional and large-scale non-convex logistic regression and robust linear regression problems. Experimental results show that ARCm significantly outperforms CR, CRm, ARC and TR, where it is 
    $10$\%-$50$\% faster than ARC (the most competitive method among CR, CRm, ARC and TR) in terms of iterations for convergence.
\end{itemize}

We really appreciate CRm, which first accelerates CR by momentum. We would like to mention our main difference with CRm \cite{wang2020cubic}. Firstly, the scheme of the momentum for ARCm is different from that in CRm due to the adaptive selection strategy for $M_{k}$. Secondly, the step size for momentum in ARCm is free of gradient evaluation but CRm requires. Thirdly, besides the global convergence, we also study the local convergence of ARCm under the \KL property, which is more general than the local error-bound condition studied for CRm. Furthermore, we analyzed the ARCm with inexact CRS solutions, which is more applicable in real applications.

\section{The Proposed ARCm Algorithm}

The detailed pseudocode for ARCm is shown in Algorithm \ref{alg_arcm}. Here, we first assume that the cubic subproblem in Step 2 of Algorithm 1 is exactly computed. In Section \ref{sec_arcm_inexact}, we will analyze the convergence property of ARCm with inexact solutions in Step 2. 

Firstly, the momentum $\bm{v}_{k-1}$ is an aggregation of previously accepted descent steps. It is always uniformly bounded if $\beta_{k} \leq \tau <1$ and $\bm{s}_{k}$ is bounded (we will prove it in Lemma \ref{lemma_4}) for all $k \leq T$. At the beginning of the algorithm (i.e., $k \ll T$), $\|\bm{s}_{k}\|_2$ is usually relative large. Therefore, $\alpha_1 \|\bm{s}_{k}\|_2$ and $\alpha_2 \|\bm{s}_{k}\|_2^2$ dominate the term $\min \left(\tau, \alpha_1 \|\bm{s}_k\|_2, \alpha_2 \|\bm{s}_k\|_2^2 \right)$, i.e., $\tau = \min \left(\tau, \alpha_1 \|\bm{s}_k\|_2, \alpha_2 \|\bm{s}_k\|_2^2 \right)$ is a popular choice for step size of momentum. When $\bm{x}_{k}$ approaches the local minima (then $\|\bm{s}_k\|_2 \approx 0$), we may not expect momentum with large step size to work for second-order methods. Therefore, we adopt $\alpha_1 \|\bm{s}_k\|_2$ and $\alpha_2 \|\bm{s}_k\|_2^2$ in selecting $\beta_{k}$ (where $\beta_{k} \approx 0$) in order to preserve the local convergence property of the second-order method. Secondly, we require that $f(\bm{z}_{k+1}) \leq f(\bm{y}_{k+1})$ since we do not hope to violate the sufficient descent of the objective in CR and ARC. Here, such $\bm{z}_{k+1}$ must exist (e.g., $\beta_{k}=0$) and we may use the bisection method to search appropriate $\beta_{k}$. The following  Theorem \ref{theorem_0} shows that $\bm{z}_{k+1}$ with nonzero $\beta_{k}$ may exist under two cases, where the momentum term helps the convergence of ARCm (i.e., $f(\bm{z}_{k+1}) < f(\bm{y}_{k+1})$ holds).


\begin{theorem}
\label{theorem_0}
Assume that the Hessian $\nabla^2f(\bm{x})$ of $f(\bm{x})$ on the line segment $[\bm{x}_{k}, \bm{x}_{k}+\bm{s}_{k}]$ is $L_{k}$-Lipschitz (i.e., for all $\bm{x}, \bm{y} \in [\bm{x}_{k}, \bm{x}_{k}+\bm{s}_{k}]$ we have $\left\|\nabla^2f(\bm{x}) - \nabla^2f(\bm{y})\right\|_2 \leq L_{k} \left\|\bm{x} - \bm{y} \right\|_2$), then in the following two cases the momentum may help the convergence, i.e., there exist small enough $\beta_{k}$ such that $f(\bm{z}_{k+1}) < f(\bm{y}_{k+1})$:
    
    (i) $L_{k} < \sigma_{k}$ and $\bm{s}_{k}^{\top}\bm{v}_{k-1}>0$;
    
    (ii) $L_{k} > \sigma_{k}$ and $\bm{s}_{k}^{\top}\bm{v}_{k-1}<0$.
\end{theorem}

\begin{proof}
Suppose that the Hessian of  $f(\bm{x})$ is $L$-Lipschitz on a ball centered at $\bm{y}_{k+1}$ with radius $\tau \|\bm{v}_{k-1}\|_2$, then we have 
\begin{equation*}
\begin{split}
    f(\bm{z}_{k+1}) - f(\bm{y}_{k+1}) \leq \beta \cdot \nabla f(\bm{y}_{k+1})^{\top}\bm{v}_{k-1} \\
    + \frac{1}{2} \beta^2 \cdot \bm{v}_{k-1}^{\top} \nabla^{2}f(\bm{x}_{k+1})\bm{v}_{k-1} + \frac{L}{6} \beta^3 \cdot \|\bm{v}_{k-1}\|_2^3,
\end{split}
\end{equation*}
by \cite[Lemma~1]{nesterov2006cubic} (we also provide the useful result in supplementary material (\ref{eq15})).
If $\nabla f(\bm{y}_{k+1})^{\top}\bm{v}_{k-1} < 0$, then there exists a small enough $\beta>0$ such that $f(\bm{z}_{k+1}) < f(\bm{y}_{k+1})$. In the remaining part, we show that in the above two cases, $\nabla f(\bm{y}_{k+1})^{\top}\bm{v}_{k-1} < 0$ may hold. Using the properties of the cubic regularization method \cite[Lemma~1~\&~(2.5)]{nesterov2006cubic}, we have
\begin{equation}
\label{momentun_decrease}
    \begin{split}
        & \nabla f(\bm{y}_{k+1})^{\top}\bm{v}_{k-1}\\
        = & \left(\nabla f(\bm{y}_{k+1}) - \nabla f(\bm{x}_{k}) \right)^{\top} \bm{v}_{k-1}  + \nabla f(\bm{x}_{k})^{\top} \bm{v}_{k-1} \\
        = & \left(\nabla f(\bm{y}_{k+1}) - \nabla f(\bm{x}_{k}) \right)^{\top} \bm{v}_{k-1} \\
        & - \left( \nabla^{2}f(\bm{x}_{k})\bm{s}_{k} + \frac{1}{2}\sigma_{k} \|\bm{s}_{k}\|_2 \cdot \bm{s}_{k} \right)^{\top}\bm{v}_{k-1} \\
        = & - \frac{1}{2}\sigma_{k} \|\bm{s}_{k}\|_2 \cdot \bm{s}_{k}^{\top}\bm{v}_{k}\\
        & + \left(\nabla f(\bm{y}_{k+1}) - \nabla f(\bm{x}_{k}) - \nabla^{2}f(\bm{x}_{k})\bm{s}_{k} \right)^{\top} \bm{v}_{k} \\
        \leq & - \frac{1}{2}\sigma_{k} \|\bm{s}_{k}\|_2 \cdot \bm{s}_{k}^{\top}\bm{v}_{k}  + \frac{1}{2}L_{k} \|\bm{s}_{k}\|_2 \cdot \left|\bm{s}_{k}^{\top}\bm{v}_{k}\right|.
    \end{split}
\end{equation}

If $L_{k} < \sigma_{k}$ (i.e., $L_k \leq \sigma_{k} - \epsilon$ for some $\epsilon > 0$) and $\bm{s}_{k}^{\top}\bm{v}_{k-1}>0$, then 
\begin{equation*}
\begin{split}
    \nabla f(\bm{y}_{k+1})^{\top}\bm{v}_{k-1} & \leq - \frac{1}{2}\sigma_{k} \|\bm{s}_{k}\|_2 \cdot \bm{s}_{k}^{\top}\bm{v}_{k}  + \frac{1}{2}L_{k}\|\bm{s}_{k}\|_2 \cdot \left|\bm{s}_{k}^{\top}\bm{v}_{k}\right| \\
    & \leq -\frac{1}{2} \epsilon \|\bm{s}_{k}\|_2 \cdot  \left|\bm{s}_{k}^{\top}\bm{v}_{k}\right| < 0.
\end{split}
\end{equation*}
In the last inequality of (\ref{momentun_decrease}), we use the inequality that 
\begin{equation*}
\begin{split}
    \left|\left(\nabla f(\bm{y}_{k+1}) - \nabla f(\bm{x}_{k}) - \nabla^{2}f(\bm{x}_{k})\bm{s}_{k} \right)^{\top} \bm{v}_{k} \right| \\
    \leq \frac{1}{2}L_{k} \|\bm{s}_{k}\|_2 \cdot \left|\bm{s}_{k}^{\top}\bm{v}_{k}\right|.
\end{split}
\end{equation*}
If $L_{k} > \sigma_{k}$, 
$\bm{s}_{k}^{\top}\bm{v}_{k-1}<0$ and 
$$
\left(\nabla f(\bm{y}_{k+1}) - \nabla f(\bm{x}_{k}) - \nabla^{2}f(\bm{x}_{k})\bm{s}_{k} \right)^{\top} \bm{v}_{k} = \frac{1}{2} \widetilde{L}_{k}\|\bm{s}_{k}\|_2 \cdot \bm{s}_{k}^{\top}\bm{v}_{k-1}
$$ 
for $\sigma_{k} \leq \widetilde{L}_{k} -\epsilon < L_k$ and $\epsilon>0$, then 
\begin{eqnarray*}
&    & \nabla f(\bm{y}_{k+1})^{\top}\bm{v}_{k-1}\\
  &  = & - \frac{1}{2}\sigma_{k} \|\bm{s}_{k}\|_2 \cdot \bm{s}_{k}^{\top}\bm{v}_{k}\\
    && + \left(\nabla f(\bm{y}_{k+1}) - \nabla f(\bm{x}_{k}) - \nabla^{2}f(\bm{x}_{k})\bm{s}_{k} \right)^{\top} \bm{v}_{k}\\
    &= & - \frac{1}{2}\sigma_{k} \|\bm{s}_{k}\|_2 \cdot \bm{s}_{k}^{\top}\bm{v}_{k}  + \frac{1}{2}\widetilde{L}_{k} \cdot \bm{s}_{k}^{\top}\bm{v}_{k-1}\\
    &\leq & \frac{1}{2} \epsilon \|\bm{s}_{k}\|_2 \cdot \bm{s}_{k}^{\top}\bm{v}_{k-1} < 0.
\end{eqnarray*}

\vspace{-8mm}\hspace{6.5cm}\end{proof}

\begin{remark}
The two potential cases may happen since ARCm adaptively select $M_{k} := \sigma_{k}$ by the creteria $\rho_{k}$ where $\sigma_{k}$ may be overestimated or underestimated. When $\sigma_{k} > L_{k}$, the step $\bm{s}_{k}$ is too conservative and the new step $\bm{x}_{k+1} - \bm{x}_{k}$ contributes to lower objective value if $\bm{v}_{k-1}$ is also a descent direction (i.e., $\bm{s}_{k}^{\top}\bm{v}_{k-1} > 0$). Conversely, if $\sigma_{k} < L_{k}$, the step $\bm{s}_{k}$ may be too aggressive, then the momentum with opposite direction (i.e., $\bm{s}_{k}^{\top}\bm{v}_{k-1} < 0$) may correct the imperfect step $\bm{s}_{k}$.
In CRm, the momentum $\bm{v}_{k-1}$ is of highly correlated to $\bm{s}_{k}$ that $\bm{s}_{k}^{\top}\bm{v}_{k-1}>0$ and $M > L_{k}$, then only the first case will happen. Therefore, the momentum scheme in ARCm is more appropriate for ARC since ARC is more ambitious than vanilla CR.  
\end{remark}

\subsection{Global Convergence}
To prove global convergence, the following mild assumptions are essential.

\begin{assumption}
\label{assump}
We have the following assumptions for the objective $f(\bm{x})$:
\begin{enumerate}
    \item $f(\bm{x})$ is globally second-order differentiable with respect to $\bm{x}$.
    
    \item $f(\bm{x})$ is bounded below, i.e., $f^{*} = \inf_{\bm{x}} f(\bm{x}) > -\infty$. 
    
    \item For the given initial guess $\bm{x}_0$, there exists a closed convex set $\mathcal{F}$ such that the level set $\mathcal{L}(\bm{x}_0):=\{\bm{x}: f(\bm{x}) \leq f(\bm{x}_0) \} \subseteq \mathcal{F}$. For all $\bm{x}, \tilde{\bm{x}} \in \mathcal{F}$, we have 
    \begin{equation}
        \begin{split}
            \left\| \nabla^2f(\bm{x}) \right\|_2 & \leq \kappa_{\rm{H}},\\
            \left\| \nabla f(\bm{x}) - \nabla f(\tilde{\bm{x}}) \right\|_2 & \leq L_{\rm{g}} \left\| \bm{x} - \tilde{\bm{x}} \right\|_2,\\
            \left\| \nabla^2 f(\bm{x}) - \nabla^2 f(\tilde{\bm{x}}) \right\|_2 & \leq L_{\rm{H}} \left\| \bm{x} - \tilde{\bm{x}} \right\|_2.
        \end{split}
    \end{equation}
\end{enumerate}
\end{assumption}

The second and the third assumption hold if $\lim_{\bm{x} \to \infty} f(\bm{x}) = +\infty$ (i.e., $f(\bm{x})$ is level bounded that $\mathcal{L}(\bm{x}_0)$ is bounded) and $f(\bm{x})$ is smooth enough, which is common in machine learning, e.g., non-negative (smooth) loss with ($\ell_2$) regularization terms. For the objective that is not bounded below, global convergence is usually hard to be achieved theoretically. We then present some useful propositions and lemmas in preparation for deriving global convergence. Some proofs are placed in the supplement.

\begin{proposition}
[{\cite[Lemma~4]{nesterov2006cubic}}, {\cite[Lemma~3.3]{cartis2011adaptive}}]
\label{prop1}
If $\bm{s}_k \in \arg \min_{\bm{s}} m_k(\bm{s})$, where $m_k(\bm{s}) := f(\bm{x}_k) + \nabla f(\bm{x}_k)^{\top} \bm{s} + \frac{1}{2} \bm{s}^{\top} \nabla^{2} f(\bm{x}_k) \bm{s} + \frac{\sigma_{k}}{6} \|\bm{s}\|_2^3$,
then  
\begin{equation}
    f(\bm{x}_{k}) - m_{k}(\bm{s}_{k}) \geq \frac{1}{12}\sigma_{k} \|\bm{s}_{k}\|_2^3.
\end{equation}
Furthermore, let $\mathcal{S}:=\{i: \rho_{i} > \eta_1\}$ denote the the set of index that successful update occurs, then for any $k \in \mathcal{S}$ we have
\begin{equation}
    f(\bm{x}_{k}) - f(\bm{y}_{k+1}) \geq \eta_1 \cdot \left( f(\bm{x}_{k}) - m_{k}(\bm{s}_{k}) \right) \geq \frac{\eta_1}{12} \sigma_{k}  \|\bm{s}_{k}\|_2^3.
\end{equation}
\end{proposition}

\begin{lemma}
\label{lemma_1}
Under Assumption \ref{assump}, the adaptive penalty parameter $\sigma_{k}$ cannot be arbitrarily large, i.e., 
\begin{equation}
    \sigma_{k} \leq \max\left\{L_{\rm{H}} \gamma_1, \sigma_{\min}\right\}:=\sigma_{\max}.
\end{equation}
\end{lemma}

\begin{lemma}
\label{lemma_3}
Denote 
    $\mathcal{S}_{j}:= \left\{ i < j: \rho_{i}> \eta_1 \right\}$,
and
    $\mathcal{U}_{j}:= \left\{ i < j: \rho_{i} \leq \eta_1 \right\}$. 
Then 
\begin{equation}
\label{eq1}
    \left|\mathcal{U}_{j}\right| \leq \left \lceil \log \frac{\max\left\{L_{\rm{H}} \gamma_1, \sigma_{\min}\right\}}{\sigma_{\min}} \right \rceil \cdot \left | \mathcal{S}_{j} \right |,
\end{equation}
and
\begin{equation}
\label{eq2}
    \left | \mathcal{S}_{T} \right | \geq \frac{T}{1+\left \lceil \log \frac{\max\left\{L_{\rm{H}} \gamma_1, \sigma_{\min}\right\}}{\sigma_{\min}} \right \rceil}.
\end{equation}
\end{lemma}

\begin{proof}
If the current update is successful, then we need at most $\left \lceil \log \frac{\sigma_{\max}}{\sigma_{\min}}\right \rceil$ steps of unsuccessful updates to achieve the next successful update, according to Lemma \ref{lemma_1}. Combining with (\ref{eq1}) and the fact that $|\mathcal{S}_{T}|+|\mathcal{U}_{T}| = T$, we have the relation (\ref{eq2}).
\end{proof}

\begin{lemma}
\label{lemma_4}
The followings hold for $\bm{s}_{k}$ and $\mathcal{S}_{T}$:
\begin{equation*}
    \max_{i \in \mathcal{S}_{T}} \|\bm{s}_{i}\|_2 \leq \left( \frac{12(f(\bm{x}_0 ) - f^{*})}{\eta_1 \sigma_{\min}} \right)^{1/3},
\end{equation*}
and
\begin{equation*}
    \min_{i \in \mathcal{S}_{T}} \|\bm{s}_{i}\|_2 \leq \left( \frac{12(f(\bm{x}_0 ) - f^{*})}{\left|\mathcal{S}_{T}\right| \eta_1 \sigma_{\min}} \right)^{1/3}.
\end{equation*}
Therefore, we have
\begin{equation*}
    \left\| \bm{v}_{k} \right\|_2 \leq \frac{1}{1-\tau}  \left( \frac{12(f(\bm{x}_0 ) - f^{*})}{\eta_1 \sigma_{\min}} \right)^{1/3},
\end{equation*}
for all $k \in \mathcal{S}_{T}$.
\end{lemma}

\begin{proof}
According to Proposition \ref{prop1}, we have 
\begin{eqnarray}
\label{eq14}
   & & \sum_{k \in \mathcal{S}_{T}}\frac{1}{12} \eta_1 \sigma_{k} \|\bm{s}_{k}\|_2^3\\
 \nonumber  & \leq & \sum_{k \in \mathcal{S}_{T}} f(\bm{x}_{k}) - f(\bm{y}_{k+1}) \leq \sum_{k \in \mathcal{S}_{T}} f(\bm{x}_{k}) - f(\bm{x}_{k+1})\\
 \nonumber  & = & \sum_{k=0}^{T} f(\bm{x}_{k}) - f(\bm{x}_{k+1}) \leq f(\bm{x}_{0}) - f^{*},
\end{eqnarray}
where the second inequality holds since $f(\bm{x}_{k+1}) \leq f(\bm{y}_{k+1})$, and the equality is satisfied because $f(\bm{x}_{k+1}) = f(\bm{x}_{k})$ for all $k \notin \mathcal{S}_{T}$. Then, we have
\begin{equation*}
    \sum_{k \in \mathcal{S}_{T}} \|\bm{s}_{k}\|_2^3 \leq \frac{12(f(\bm{x}_0 ) - f^{*})}{\eta_1 \sigma_{\min}}
\end{equation*}
as $\sigma_{k} \geq \sigma_{\min}$. Therefore, the first two inequalities hold by using $\max_{i \in \mathcal{S}_{T}} \|\bm{s}_{i}\|_2^3 \leq \sum_{i \in \mathcal{S}_{T}} \|\bm{s}_{i}\|_2^3$ and $\min_{i \in \mathcal{S}_{T}} \|\bm{s}_{i}\|_2^3 \leq \frac{1}{|\mathcal{S}_{T}|} \sum_{i \in \mathcal{S}_{T}} \|\bm{s}_{i}\|_2^3$. The last inequality in the lemma holds since for any $k \in \mathcal{S}_{T}$, 
\begin{equation*}
    \|\bm{v}_{k}\|_2 \leq \frac{1}{1-\tau} \cdot \max_{i\leq k, i \in \mathcal{S}_{T}} \|\bm{s}_{i}\|_2 \leq \frac{1}{1-\tau} \cdot \max_{i \in \mathcal{S}_{T}} \|\bm{s}_{i}\|_2.
\end{equation*}

\vspace{-8mm}\hspace{7cm}\end{proof}

\begin{lemma}
\label{lemma_2}
If $k \in \mathcal{S}_{T}$, we have
\begin{equation}
    \left\| \nabla f(\bm{x}_{k+1}) \right\|_2 \leq c_1 \left\| \bm{s}_{k} \right\|_2^2,
\end{equation}
and
\begin{equation}
    \lambda_{\min}\left(\nabla^2 f(\bm{x}_{k+1})\right) \geq -c_2 \left\| \bm{s}_{k} \right\|_2,
\end{equation}
where $c_1 = \frac{1}{2}\max\left\{L_{\rm{H}} \gamma_1, \sigma_{\min}\right\} + \frac{1}{2} L_{\rm{H}} + \frac{\alpha_2 L_{\rm{g}}}{1-\tau} \left( \frac{12(f(\bm{x}_0 ) - f^{*})}{\eta_1 \sigma_{\min}} \right)^{1/3}$ and $c_2 = \frac{1}{2}\max\left\{L_{\rm{H}} \gamma_1, \sigma_{\min}\right\} + L_{\rm{H}} + \frac{\alpha_1 L_{\rm{H}}}{1-\tau} \left( \frac{12(f(\bm{x}_0 ) - f^{*})}{\eta_1 \sigma_{\min}} \right)^{1/3}$. 
\end{lemma}

\begin{proof}
If $k \in \mathcal{S}_{T}$, we have 
\begin{equation*}
    \left\| \nabla f(\bm{y}_{k+1}) \right\|_2 \leq \frac{1}{2}\left(\sigma_{k}+L_{\rm{H}}\right) \left\| \bm{s}_{k} \right\|_2^2,
\end{equation*}
and
\begin{equation*}
    \lambda_{\min}\left(\nabla^2 f(\bm{y}_{k+1})\right) \geq -\left( \frac{1}{2}\sigma_{k} + L_{\rm{H}} \right) \left\| \bm{s}_{k} \right\|_2.
\end{equation*}
The proof can be found in \cite[Lemma~3~\&~5]{nesterov2006cubic}. 
We then derive the error bounds for $\left\| \nabla f(\bm{x}_{k+1}) \right\|_2$ and $\lambda_{\min}(\nabla^2 f(\bm{x}_{k+1}))$. 
\begin{equation*}
    \begin{split}
        & \left \| \nabla f(\bm{x}_{k+1}) \right \|_2\\
        \leq & \left \| \nabla f(\bm{y}_{k+1}) \right \|_2 + \left \| \nabla f(\bm{x}_{k+1}) - \nabla f(\bm{y}_{k+1}) \right \|_2\\
        \leq & \left \| \nabla f(\bm{y}_{k+1}) \right \|_2 + L_{\rm{g}} \beta_{k} \|\bm{v}_{k}\|_2\\
        \leq & \frac{1}{2}\left(\sigma_{k}+L_{\rm{H}}\right) \left\| \bm{s}_{k} \right\|_2^2 + L_{\rm{g}} \alpha_2 \|\bm{s}_{k}\|_2^2 \cdot \|\bm{v}_{k}\|_2 
        \leq  c_1 \|\bm{s}_{k}\|_2^2,
    \end{split}
\end{equation*}
where $c_1 = \frac{1}{2}\sigma_{\max} + \frac{1}{2} L_{\rm{H}} + \frac{\alpha_2 L_{\rm{g}}}{1-\tau} \left( \frac{12(f(\bm{x}_0 ) - f^{*})}{\eta_1 \sigma_{\min}} \right)^{1/3}$. 
Furthermore, we have
\begin{equation*}
    \begin{split}
        & \lambda_{\min}\left(\nabla^2 f(\bm{x}_{k+1})\right)\\
        \geq & \lambda_{\min}\left(\nabla^2 f(\bm{y}_{k+1})\right) - \left\| \nabla^2 f(\bm{x}_{k+1}) - \nabla^2 f(\bm{y}_{k+1}) \right\|_2\\
        \geq & \lambda_{\min}\left(\nabla^2 f(\bm{y}_{k+1})\right) - L_{\rm{H}} \left\| \bm{x}_{k+1} - \bm{y}_{k+1} \right\|_2\\
        \geq & \lambda_{\min}\left(\nabla^2 f(\bm{y}_{k+1})\right) - L_{\rm{H}} \beta_{k} \|\bm{v}_{k}\|_2 \\
        \geq & -\left( \frac{1}{2}\sigma_{k} + L_{\rm{H}} \right) \left\| \bm{s}_{k} \right\|_2 - L_{\rm{H}} \alpha_1 \|\bm{s}_{k}\|_2 \cdot \|\bm{v}_{k}\|_2\\
        \geq & - c_2 \|\bm{s}_{k}\|_2,
    \end{split}
\end{equation*}
where $c_2 = \frac{1}{2}\sigma_{\max} + L_{\rm{H}} + \frac{\alpha_1 L_{\rm{H}}}{1-\tau} \left( \frac{12(f(\bm{x}_0 ) - f^{*})}{\eta_1 \sigma_{\min}} \right)^{1/3}$. 
\end{proof}

\begin{theorem}
\label{theorem_1}
We introduce the following measure of the local optimality:
\begin{equation}
    \mu(\bm{x}) = \max \left \{\sqrt{\frac{1}{c_1} \left\|\nabla f(\bm{x})\right\|}, -\frac{1}{c_2}  \lambda_{\min}\left(\nabla^2 f(\bm{x})\right) \right\},
\end{equation}
where $c_1>0$ and $c_2>0$ are two universal constant defined in Lemma \ref{lemma_2}.
Under Assumption \ref{assump}, let the sequence $\{\bm{x}_{k}\}_{k=1}^{T}$ be generated by Algorithm \ref{alg_arcm}, then
\begin{equation}
\begin{split}
    & \min_{1 \leq k \leq T} \mu(\bm{x}_{k}) \\
    & \leq \left( \frac{12(f(\bm{x}_0 ) - f^{*})}{\eta_1 \sigma_{\min} T} \left (1+\left \lceil \log \frac{\max\left\{L_{\rm{H}} \gamma_1, \sigma_{\min}\right\}}{\sigma_{\min}} \right \rceil \right) \right)^{1/3}\\
    & = \mathcal{O}\left(T^{-1/3} \right).
\end{split}
\end{equation}
\end{theorem}

\begin{proof}
Lemma \ref{lemma_2} implies that $\mu(\bm{x}_{k+1}) \leq \|\bm{s}_{k}\|_2$ for all $k \in \mathcal{S}_{T}$. Then we have $\min_{1 \leq k \leq T} \mu(\bm{x}_{k}) \leq \min_{k \in \mathcal{S}_{T}} \|\bm{s}_{k}\|_2$.  Combining it with Lemma \ref{lemma_3} and Lemma \ref{lemma_4}, we finish the proof. 
\end{proof}

\begin{remark}
Note that the proposed ARCm also satisfies the global convergence rate $\mathcal{O}(T^{-1/3})$ for the measure $\mu(\cdot)$, as the vanilla CR \cite{nesterov2006cubic} and ARC \cite{cartis2011adaptive}. We do not expect to improve the convergence rate since ARCm is still a second-order method.
\end{remark}

\subsection{Local Convergence}
In this subsection, we let $T = +\infty$.  If the accumulation point $\bar{\bm{x}}$ of the sequence $\{\bm{x}_{k}\}_{k=0}^{\infty}$ generated by CR satisfies $\lambda_{\min}(\nabla^2 f(\bar{\bm{x}}))>0$, then it achieves the local quadratic convergence \cite{nesterov2006cubic}. However, the Hessian $\nabla^2 f(\bm{x})$ at the local optima is usually not necessary to be positive definite. Therefore, the local quadratic convergence may not work. Yue et al. \cite{yue2019quadratic} proved that CR achieves the quadratic convergence if the objective $f(\bm{x})$ satisfies the local error bound condition, which is weaker than the local positive definiteness of Hessian. Later, Zhou et al. \cite{zhou2018convergence} generalized the results in \cite{yue2019quadratic} to the objective that satisfies the \KL property. Here, we show that the proposed ARCm also achieves the local convergence under the \KL property. In the following analysis, sets $\mathcal{S}:=\{i: \rho_{i} > \eta_1\}$ and $\mathcal{U}:=\{i: \rho_{i} \leq \eta_1\}$ are in the ascending order. Lemma \ref{lemma_3} implies that $|\mathcal{S}|$ must be infinite but $|\mathcal{U}|$ may be finite.  We use $k_j \in \mathcal{S}$ to denote the $j$-th element of $\mathcal{S}$. 

\begin{definition}
A differentiable function $f(\cdot)$ is said to satisfy the \KL property if for any compact set $\bar{\mathcal{X}}$ where $f(\cdot)$ takes a constant value $\bar{f}$, there exist $\epsilon_1, \epsilon_2 >0$ such that for all $\bar{\bm{x}} \in \bar{\mathcal{X}}$ and $\bm{x} \in \{\bm{z}: \text{dist}(\bm{z}, \bar{\mathcal{X}})<\epsilon_1, \bar{f} < f(\bm{z})<\bar{f}+\epsilon_2\}$,
\begin{equation}
\label{eq3}
    \phi^{\prime}(f(\bm{x}) - \bar{f})\left\| \nabla f(\bm{x}) \right\| \geq 1 
\end{equation}
holds, where $\phi(t) = \frac{c}{\theta}t^{\theta}$ for some $c>0$ and $\theta \in (0,1)$. Then, (\ref{eq3}) is equivalent to
\begin{equation}
\label{eq4}
    f(\bm{x}) - \bar{f} \leq c_0\left\| \nabla f(\bm{x}) \right\|^{\frac{1}{1-\theta}}
\end{equation}
with $c_0 = c^{1/(1-\theta)}$.
\end{definition}

Besides Assumption \ref{assump}, we need further but mild assumptions for the local convergence. As is discussed following Assumption \ref{assump}, the level boundedness (Assumption \ref{assump2}) and smoothness (e.g., $f(\bm{x})$ is third-order differentiable with $\bm{x}$) imply the second and the third assumptions in Assumption \ref{assump}.

\begin{assumption}
\label{assump2}
We further assume that $f(\bm{x})$ is level bounded, i.e., the level set $\mathcal{L}(\tilde{\bm{x}})$ is bounded, $\forall \tilde{\bm{x}} \in \mathcal{F}$. 
\end{assumption}

Before deriving the local convergence for ARCm under the \KL property, we first show that the sequence generated by $\{\bm{x}_{k}\}_{k}$ is Cauchy and convergent to a second-order critical point. The proofs for Lemma \ref{lemma_8} and Theorem \ref{theorem_2} are extended from \cite{zhou2018convergence} and we put the tedious details in the supplementary material due to the page limit.  

\begin{lemma}
\label{lemma_8}
Under Assumption \ref{assump} and \ref{assump2}, the followings hold for the sequence $\{\bm{x}_{k}\}_{k=0}^{+\infty}$ generated by ARCm (Algorithm \ref{alg_arcm}):
\begin{enumerate}
    \item $\bar{f} := \lim_{k \to \infty} f(\bm{x}_{k})$ exists. 
    
    \item The sequence $\{\bm{x}_{k}\}_{k=0}^{+\infty}$ is bounded and $\lim_{k \to \infty}\|\bm{x}_{k+1} - \bm{x}_{k}\|_2 = 0$. Moreover, the set $\bar{\mathcal{X}}$ of accumulation points of the sequence is non-empty, satisfying
    \begin{equation*}
        f(\bar{\bm{x}}) = \bar{f}, \quad \nabla f(\bar{\bm{x}}) = \bm{0}, \quad \lambda_{\min}\left(\nabla^2 f(\bar{\bm{x}})\right) \geq \bm{0},
    \end{equation*}
    for all $\bar{\bm{x}} \in \bar{\mathcal{X}}$. 
    
    \item If $f(\bm{x})$ satisfies the \KL property, then $\bar{\mathcal{X}}=\{\bar{\bm{x}}\}$ is a singleton.
\end{enumerate}
\end{lemma}

\begin{theorem}
\label{theorem_2}
Let the objective $f(\bm{x})$ satisfies Assumption \ref{assump}, \ref{assump2} and the \KL property, then there exists a large enough $j_0 \in \mathbb{N}$ such that the sequence $\{\bm{x}_{k}\}_{k=0}^{+\infty}$ (or $\{\bm{x}_{k_j+1}\}_{j=0}^{+\infty}$, $k_j \in \mathcal{S}$) generated by ARCm (Algorithm \ref{alg_arcm}) satisfies
\begin{enumerate}
    \item If $\theta \in (\frac{1}{3},1)$, then the local convergence is super-linear with
    \begin{equation*}
        \|\bm{x}_{k_{j}+1} - \bar{\bm{x}}\|_2 \leq \mathcal{O}\left(\exp\left(-\left( \frac{2\theta}{1-\theta}\right)^{j-j_0}\right) \right)
    \end{equation*}
    and
    \begin{equation*}
        \|\bm{x}_{k+1} - \bar{\bm{x}}\|_2 \leq \mathcal{O}\left(\exp\left(-\left( \frac{2\theta}{1-\theta}\right)^{\left\lceil\frac{k-k_{j_0}}{c_{3}}\right\rceil}\right) \right),
    \end{equation*}
    where $c_{3} = 1+\left \lceil \log \frac{\max\left\{L_{\rm{H}} \gamma_1, \sigma_{\min}\right\}}{\sigma_{\min}} \right \rceil$.
    
    \item If $\theta =\frac{1}{3}$, then the local convergence is linear with
    \begin{equation*}
        \|\bm{x}_{k_j+1} - \bar{\bm{x}}\|_2 \leq  \mathcal{O}\left(\exp \left(-c_{4}(j-j_0) \right) \right)
    \end{equation*}
    and
    \begin{equation*}
        \|\bm{x}_{k+1} - \bar{\bm{x}}\|_2 \leq \mathcal{O}\left(\exp \left(-c_{4}\left \lceil\frac{k-k_{j_0}}{c_{3}} \right\rceil \right) \right),
    \end{equation*}
    for some constants $c_{4}>0$.

    \item If $\theta \in (0, \frac{1}{3})$, then the local convergence is sub-linear with
    \begin{equation*}
        \|\bm{x}_{k_j+1} - \bar{\bm{x}}\|_2 \leq \mathcal{O}\left(\left(j-j_0\right)^{-\frac{2\theta}{1-3\theta}}\right)
    \end{equation*}
    and
    \begin{equation*}
        \|\bm{x}_{k+1} - \bar{\bm{x}}\|_2 \leq \mathcal{O}\left( \left\lceil\frac{k-k_{j_0}}{c_{3}}\right\rceil^{-\frac{2\theta}{1-3\theta}}\right).
    \end{equation*}
    
\end{enumerate}
\end{theorem}

\section{ARCm with Inexact Solutions}
\label{sec_arcm_inexact}

A popular approach to the exact solution of CRS in 
Step 2 of Algorithm 1, is solving the corresponding secular equation \cite{nesterov2006cubic}, where full eigendecomposition for the Hessian matrix is required (with computational complexity $\mathcal{O}(d^{3})$). Cartis et al. \cite{cartis2011adaptive} designed a Newton-Cholesky iteration method for solving CRS, however, its computational cost is still of   $\mathcal{O}(d^{3})$. 
Exact solutions of CRS are usually very computationally expensive for high-dimensional problems. In practice, inexact solvers are more popular. Cartis et al. \cite{cartis2011adaptive} projected the CRS to a Krylov subspace in order to lower the dimension. The convergence of the Krylov subspace method for inexactly solving CRS was analyzed by Carmon and Duchi \cite{carmon2018analysis}. Later, they showed that the simple gradient descent with proper learning rates finds the global solution of CRS \cite{carmon2019gradient}. Jiang et al. \cite{jiang2021accelerated} proposed an accelerated first-order method that reformulates the CRS into a constrained convex problem.  Recently, Gao et al. \cite{gao2022approximate} suggested solving an approximate secular equation rather than the exact secular equation, where partial eigendecomposition is required. Suppose that we solve the CRS by inexact solvers that satisfy Condition \ref{condition1}: 
\begin{equation*}
    \begin{split}
        & \tilde{\bm{s}}_{k} \approx \arg \min_{\bm{s}} m_{k}(\bm{s}),\\
        & \bm{y}_{k+1} = \bm{x}_{k+1} + \tilde{\bm{s}}_{k},\\
        & \tilde{\bm{v}}_{k} = \tilde{\beta}_{k}\tilde{\bm{v}}_{k-1} +  \tilde{\bm{s}}_{k}~\text{and}~\bm{z}_{k+1} = \bm{x}_{k} + \tilde{\bm{v}}_{k}\\ 
        & \hspace{2em}~\text{with}~\tilde{\beta}_{k} \in [0,\min\{\tau, \alpha_1\|\tilde{\bm{s}}_{k}\|_2, \alpha_2 \|\tilde{\bm{s}}_{k}\|_2^2\}]\\
        & \hspace{2em}~\text{and}~ f(\bm{z}_{k+1}) \leq f(\bm{y}_{k+1}),\\
        & \bm{x}_{k+1} = \bm{z}_{k+1}.
    \end{split}
\end{equation*}

\begin{condition}
\label{condition1}
Suppose that $\tilde{\bm{s}}_{k}$ is the inexact solution of the CRS in the $k$-th iteration, satisfying the following $\delta_k$-conditions:
\begin{enumerate}
    \item $m_{k}(\tilde{\bm{s}}_{k}) - f(\bm{x}_{k}) \leq -\frac{\sigma_k}{12}\|\tilde{\bm{s}}_{k}\|_2^3 + \delta_{k} < 0$;
    
    \item $\nabla m_{k}(\tilde{\bm{s}}_{k}) \leq \delta_{k}^{2/3}$;
    
    \item $\left| \left\|\tilde{\bm{s}}_{k}\right\|- \left\|\bm{s}_{k}\right\| \right| \leq \delta_{k}^{1/3}$.
\end{enumerate}
\end{condition}

Note that the first item $m_{k}(\tilde{\bm{s}}_{k}) - f(\bm{x}_{k}) < 0$ is easy to be satisfied since the Cauchy point method guaranteed the sufficient decrease except at a stable point. Before analyzing the global convergence of ARCm with inexact CRS solutions, we first provide some useful lemmas where some proofs are in the supplement.

\begin{lemma}
\label{lemma_7}
Under Assumption \ref{assump} and Condition \ref{condition1}, the adaptive penalty parameter $\sigma_{k}$ cannot be arbitrary large, i.e., 
\begin{equation*}
    \sigma_{k} \leq \max \left\{L_{\rm{H}}\gamma_1, \sigma_{\min} \right\} = \sigma_{\max}.
\end{equation*}
Therefore, (\ref{eq1}) and (\ref{eq2}) still hold for $|\mathcal{U}_{T}|$ and $|\mathcal{S}_{T}|$.
\end{lemma}

\begin{lemma}
\label{lemma_6}
Without the loss of generality, we assume that $\delta_{k} < f(\bm{x}_{0}) - f^{*}$. Suppose that Assumption \ref{assump} and Condition \ref{condition1} hold, then we have \begin{equation*}
    \max_{k \in \mathcal{S}_{T}} \|\tilde{\bm{s}}_{k}\|_2 \leq \left( \frac{24(f(\bm{x}_0 ) - f^{*})}{\eta_1 \sigma_{\min}} \right)^{1/3},
\end{equation*}
and
\begin{equation*}
    \min_{k \in \mathcal{S}_{T}} \|\tilde{\bm{s}}_{k}\|_2^3 - \frac{12\delta_{k}}{\sigma_{k}} \leq  \frac{12(f(\bm{x}_0 ) - f^{*})}{\left|\mathcal{S}_{T}\right| \eta_1 \sigma_{\min}}.
\end{equation*}
Therefore, we have
\begin{equation*}
    \left\| \tilde{\bm{v}}_{k} \right\|_2 \leq \frac{1}{1-\tau}  \left( \frac{24(f(\bm{x}_0 ) - f^{*})}{\eta_1 \sigma_{\min}} \right)^{1/3},
\end{equation*}
for all $k \in \mathcal{S}_{T}$
\end{lemma}

\begin{lemma}
\label{lemma_5}
Under Assumption \ref{assump} and Condition \ref{condition1}, if $k \in \mathcal{S}_{T}$, we have
\begin{equation}
    \left\| \nabla f(\bm{x}_{k+1}) \right\|_2 \leq c_{5} \left\| \tilde{\bm{s}}_{k} \right\|_2^2 + \left\| \nabla m_{k}(\tilde{\bm{s}}_{k}) \right\|_2,
\end{equation}
and
\begin{equation}
    \lambda_{\min}\left(\nabla^2 f(\bm{x}_{k+1})\right) \geq -c_{6} \left\| \tilde{\bm{s}}_{k} \right\|_2  - \frac{\sigma_{\max}}{2} \left|\|\bm{s}_{k}\|_2 - \|\tilde{\bm{s}}_{k}\|_2\right|,
\end{equation}
where $c_{5}=\frac{1}{2}\max\left\{L_{\rm{H}} \gamma_1, \sigma_{\min}\right\} + \frac{1}{2} L_{\rm{H}} + \frac{\alpha_2 L_{\rm{g}}}{1-\tau} \left( \frac{24(f(\bm{x}_0 ) - f^{*})}{\eta_1 \sigma_{\min}} \right)^{1/3}$ and $c_{6}=\frac{1}{2}\sigma_{\max} + L_{\rm{H}} + \frac{\alpha_1 L_{\rm{H}}}{1-\tau} \left( \frac{24(f(\bm{x}_0 ) - f^{*})}{\eta_1 \sigma_{\min}} \right)^{1/3}$
\end{lemma}

\begin{proof}
We first derive the error bound for $\left\|\nabla f(\bm{y}_{k+1}) \right\|$ and $\lambda_{\min}\left(\nabla^2 f(\bm{y}_{k+1}) \right)$: 
\begin{eqnarray}
\mbox{\hspace{-5mm}} & & \\
\nonumber  \mbox{\hspace{-5mm}} &     & \left\| \nabla f(\bm{y}_{k+1}) \right\|_2\\
   \nonumber \mbox{\hspace{-5mm}} &   \leq & \left\| \nabla f(\bm{y}_{k+1}) - \nabla f(\bm{x}_{k}) - \nabla^2 f(\bm{x}_{k}) \tilde{\bm{s}}_{k} \right\|_2\\
   \nonumber  \mbox{\hspace{-5mm}} &   & + \left\| \nabla f(\bm{x}_{k}) + \nabla^2 f(\bm{x}_{k}) \tilde{\bm{s}}_{k} + \frac{\sigma_{k}}{2}\|\tilde{\bm{s}}_{k}\|_2 \tilde{\bm{s}}_{k} \right\|_2 + \frac{\sigma_{k}}{2}\|\tilde{\bm{s}}_{k}\|_2^2\\
    \nonumber \mbox{\hspace{-5mm}} &  \leq & \frac{L_{\rm{H}}}{2}\|\tilde{\bm{s}}_{k}\|_2^2 + \left\| \nabla m_{k}(\tilde{\bm{s}}_{k}) \right\|_2 + \frac{\sigma_{k}}{2}\|\tilde{\bm{s}}_{k}\|_2^2\\
 \nonumber  \mbox{\hspace{-5mm}}   &  \leq & \left(\frac{\sigma_{\max}}{2} +\frac{L_{\rm{H}}}{2}   \right)\|\tilde{\bm{s}}_{k}\|_2^2 + \left\| \nabla m_{k}(\tilde{\bm{s}}_{k}) \right\|_2,
\end{eqnarray}
where the second inequality is due to (\ref{eq13}); and
\begin{equation*}
\begin{split}
    & \lambda_{\min}\left(\nabla^2 f(\bm{y}_{k+1}) \right)\\
    \geq & \lambda_{\min}\left(\nabla^2 f(\bm{x}_{k}) \right) - \left\|\nabla^2 f(\bm{y}_{k+1}) - \nabla^2 f(\bm{x}_{k}) \right\|_2\\
    \geq & -\frac{\sigma_{k}}{2}\|\bm{s}_{k}\|_2 - L_{\rm{H}} \|\tilde{\bm{s}}_{k}\|_2\\
    \geq & -\frac{\sigma_k}{2} \left|\|\bm{s}_{k}\|_2 - \|\tilde{\bm{s}}_{k}\|_2\right| - \frac{\sigma_k}{2} \|\tilde{\bm{s}}_{k}\|_2 - L_{\rm{H}} \|\tilde{\bm{s}}_{k}\|_2\\
    = & -\left(\frac{\sigma_{\max}}{2} - L_{\rm{H}} \right) \left\|\tilde{\bm{s}}_{k}\right\|_2 - \frac{\sigma_{\max}}{2} \left|\|\bm{s}_{k}\|_2 - \|\tilde{\bm{s}}_{k}\|_2\right|.
\end{split}
\end{equation*}
where the second inequality comes from a well-known result of cubic regularization \cite[Proposition~1]{nesterov2006cubic}. Using similar arguments as in Lemma \ref{lemma_2}, we complete the proof. We put the remaining details in the supplement.
\end{proof}

\begin{figure*}[t!]
  \centering
    \subfloat[ARCENE]{
    \label{logistic_grad_arcene}
    \includegraphics[scale=0.3]{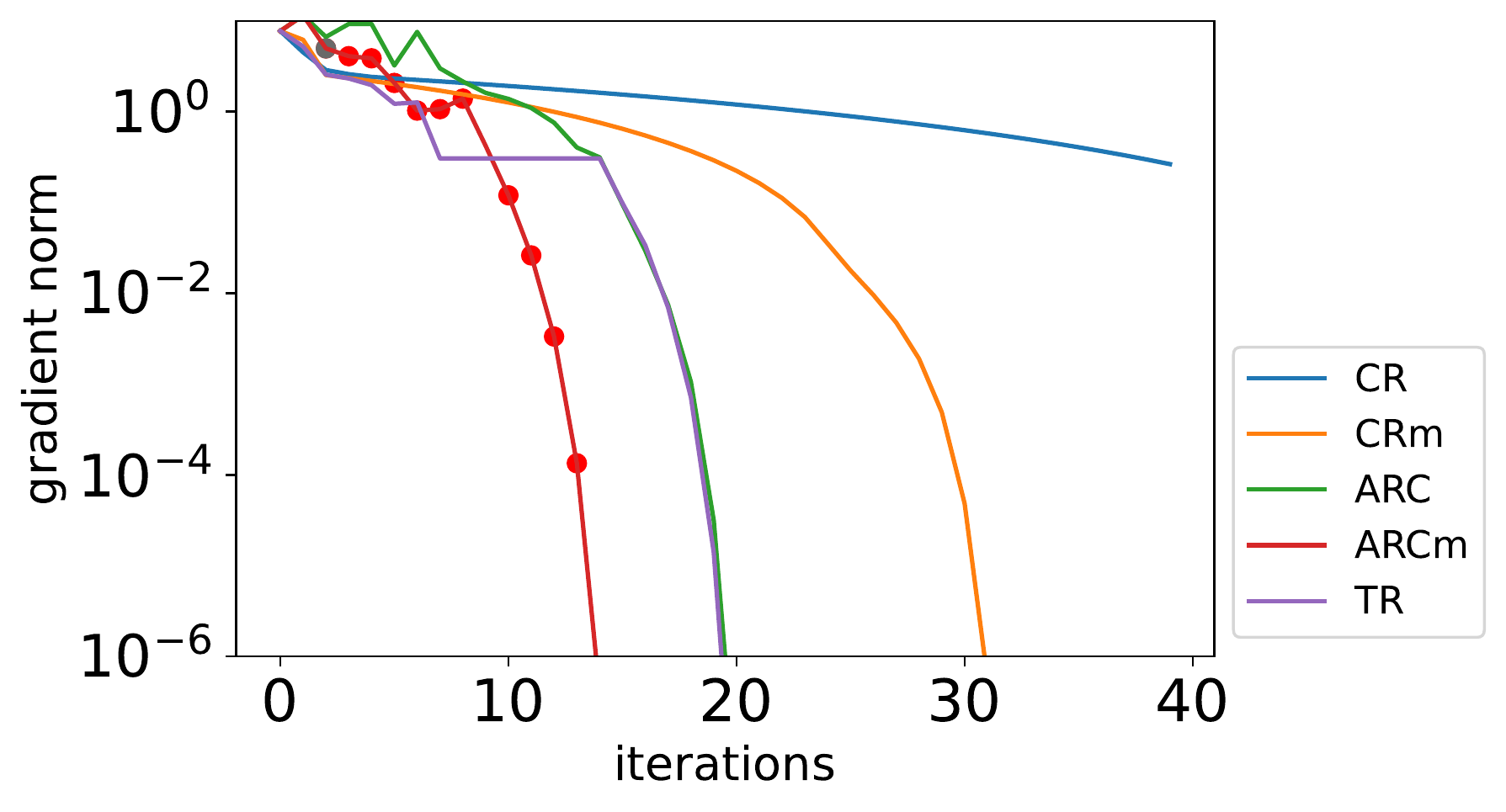}
    }\
    \subfloat[Covtype]{
    \label{logistic_grad_covtype}
    \includegraphics[scale=0.3]{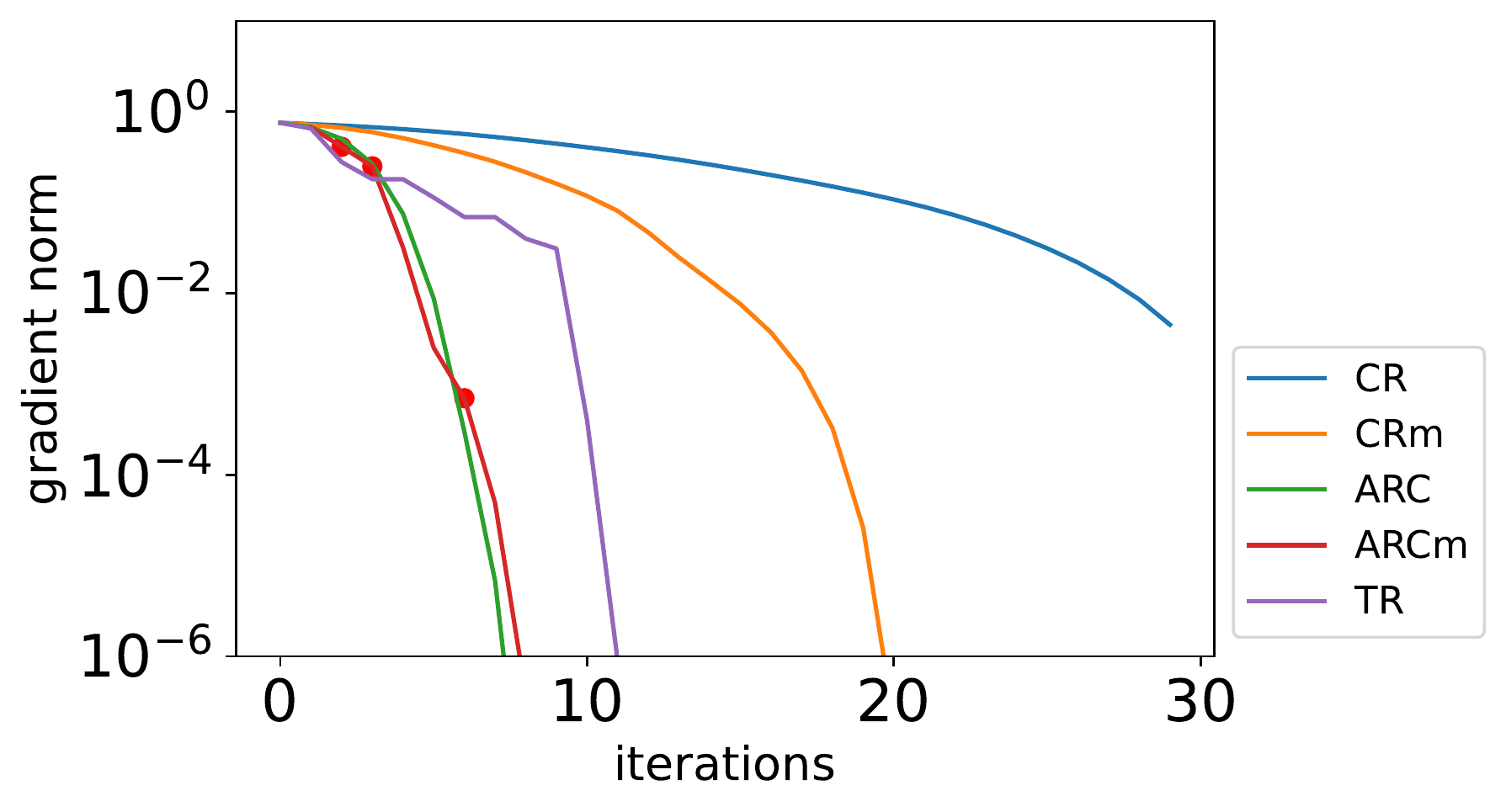}
    }\
    \subfloat[DrivFace]{
    \label{logistic_grad_drivface}
    \includegraphics[scale=0.3]{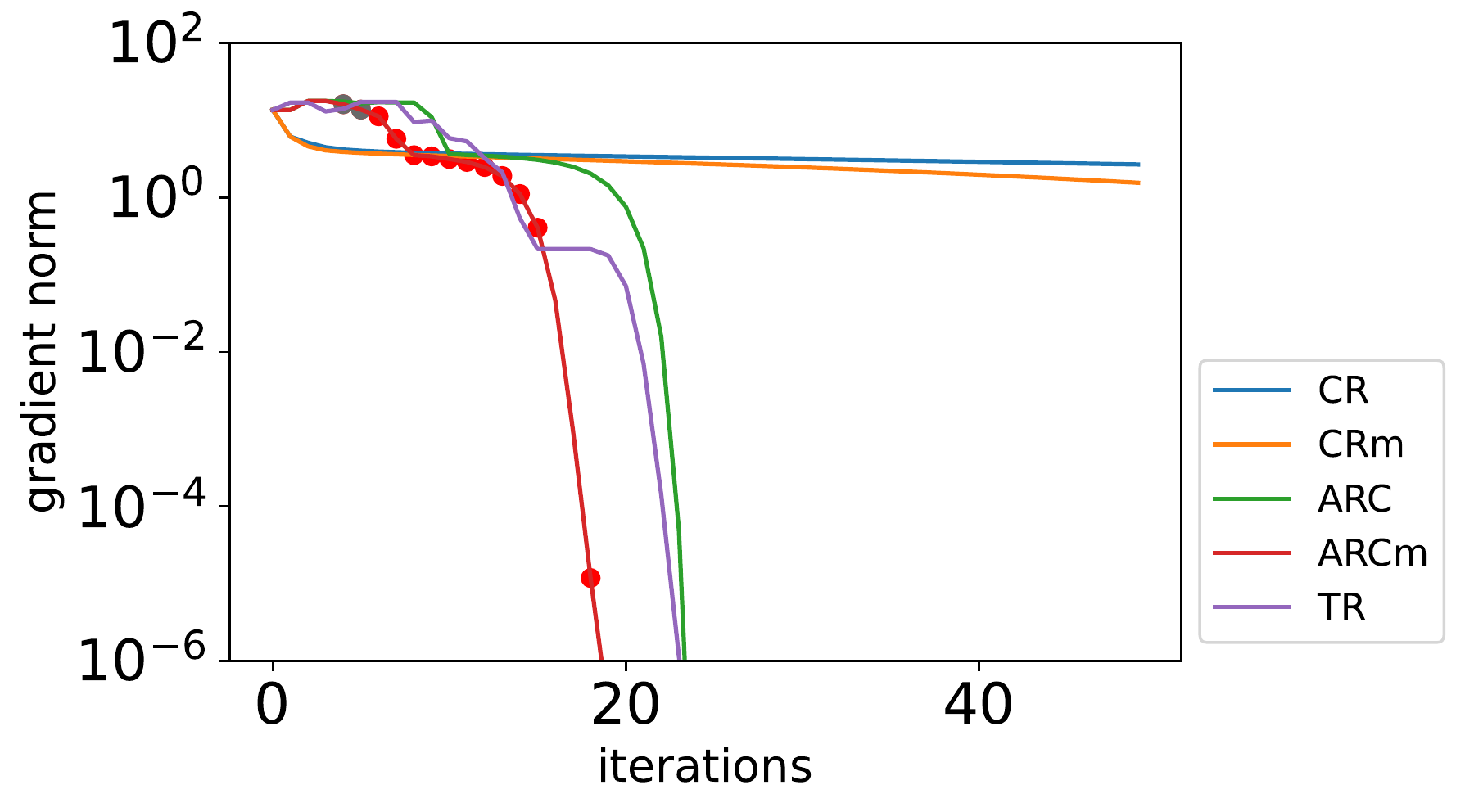}
    }\\
    \subfloat[ARCENE]{
    \label{robust_grad_arcene}
    \includegraphics[scale=0.3]{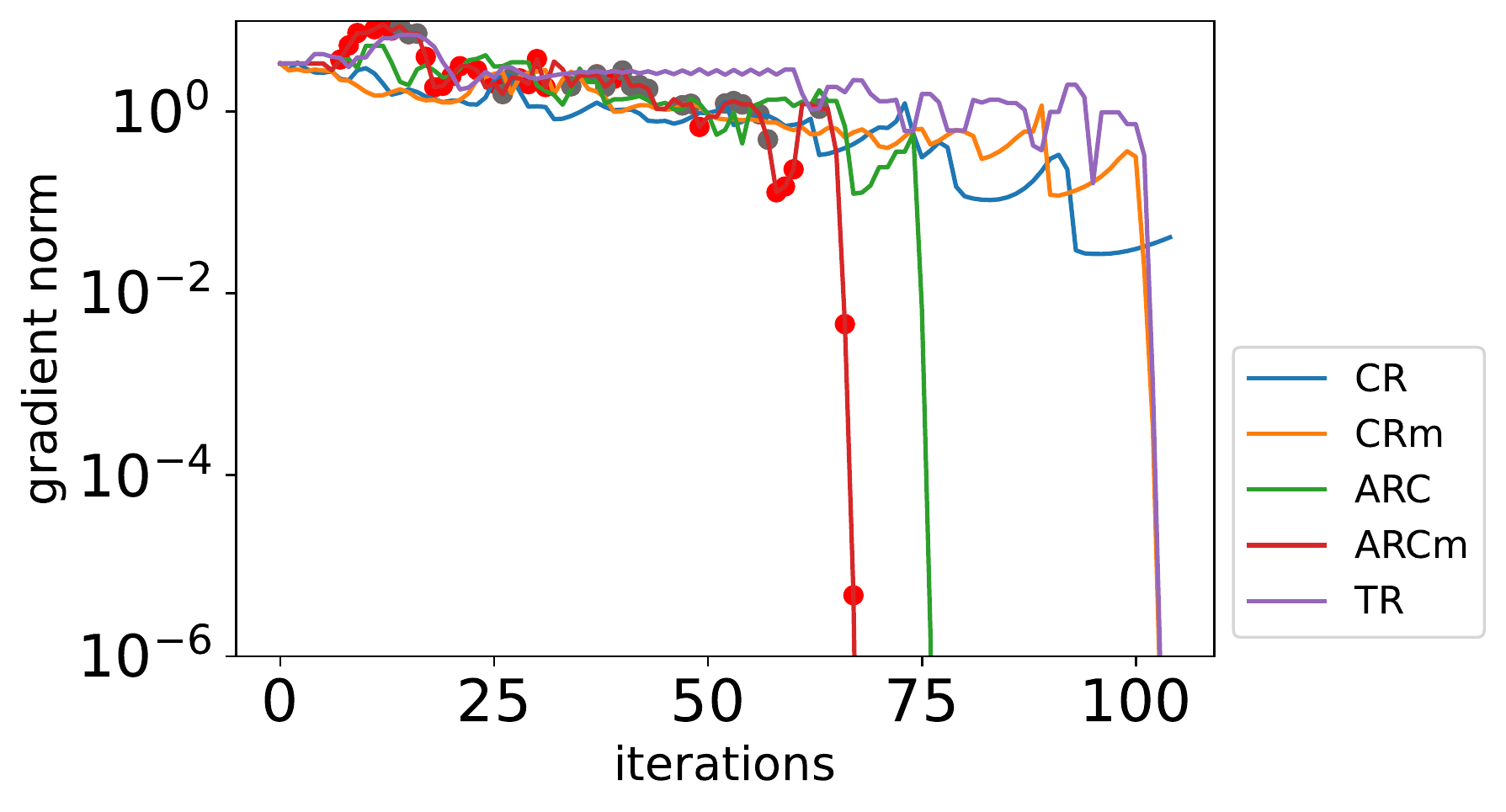}
    }\
    \subfloat[Covtype]{
    \label{robust_grad_covtype}
    \includegraphics[scale=0.3]{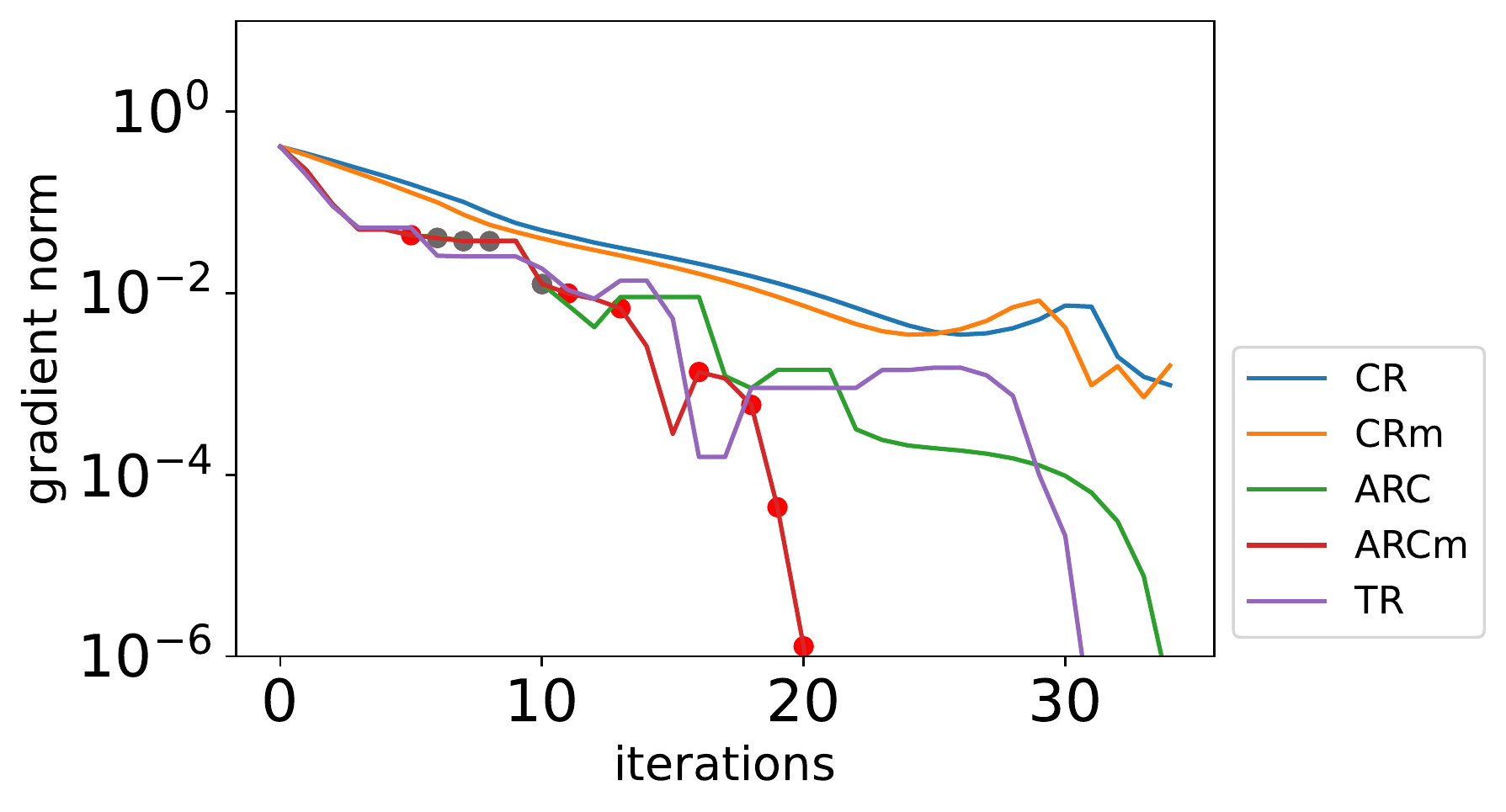}
    }\
    \subfloat[DrivFace]{
    \label{robust_grad_drivface}
    \includegraphics[scale=0.3]{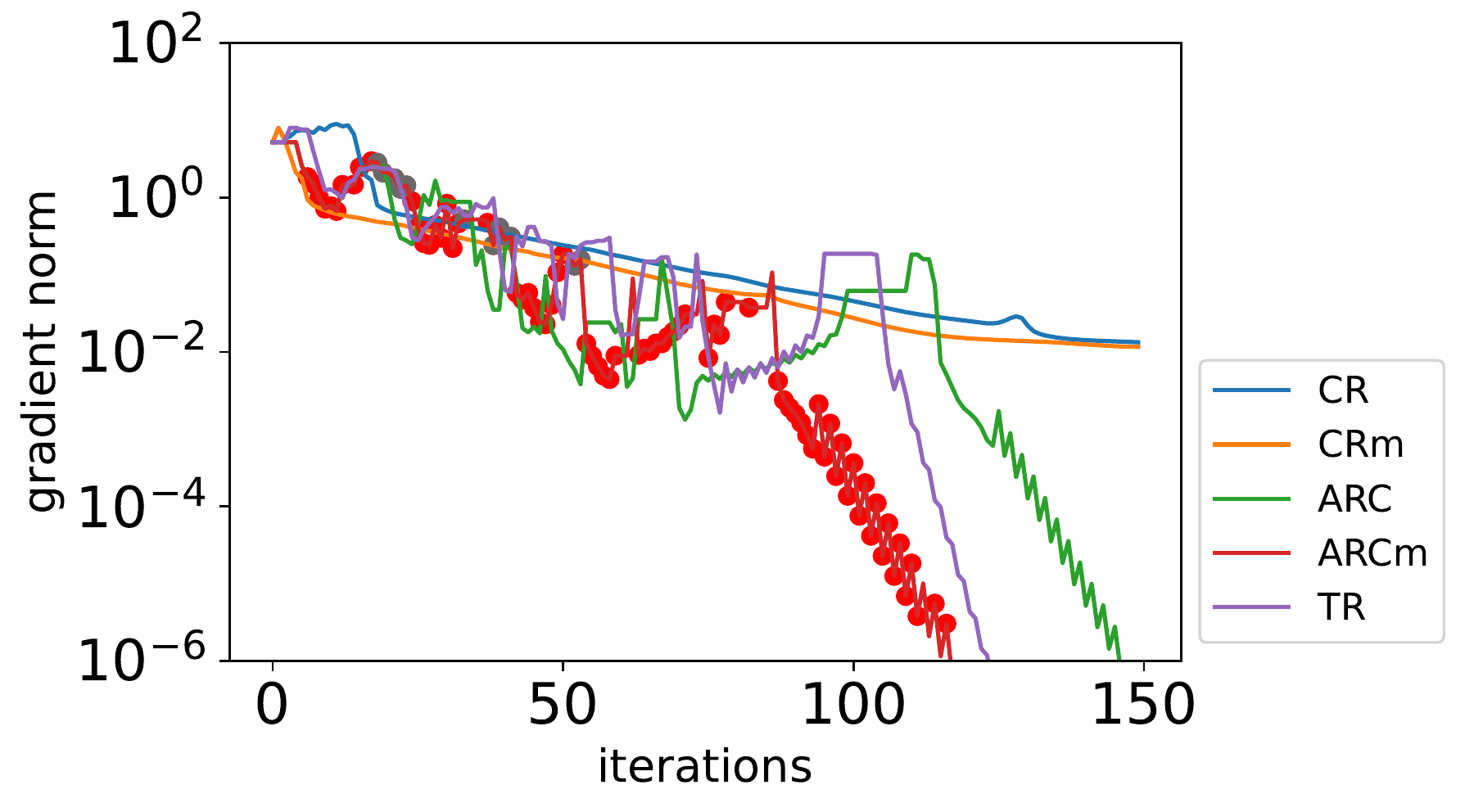}
    }\
    \caption{Gradient norm versus the iteration steps on logistic regression (the first row) and robust linear regression (the second row) models for three datasets. Red dots: $\beta_{k}>0$ and $\bm{s}_{k}^{\top}\bm{v}_{k-1} >0$; gray dots: $\beta_{k}>0$ and $\bm{s}_{k}^{\top}\bm{v}_{k-1} < 0$.}
    \label{fig_grad}
\end{figure*}

\begin{theorem}
\label{theorem_3}
We introduce the following measure of the local optimality:
\begin{equation}
    \tilde{\mu}(\bm{x}) = \max \left \{\sqrt{\frac{1}{c_{5}} \left\|\nabla f(\bm{x})\right\|}, -\frac{1}{c_{6}}  \lambda_{\min}\left(\nabla^2 f(\bm{x})\right) \right\}.
\end{equation}
Under Assumption \ref{assump} and Condition \ref{condition1}, let the sequence $\{\bm{x}_{k}\}_{k=1}^{T}$ be generated by Algorithm \ref{alg_arcm}, we have
\begin{enumerate}
    \item If $0< \delta_{k} \leq \delta$ for some $\delta>0$, then
    \begin{equation}
        \min_{1 \leq k \leq T} \tilde{\mu}(\bm{x}_{k}) \leq \mathcal{O}\left(T^{-1/3}+ \delta^{1/3} \right).
    \end{equation}
    
    \item If there exists $0<\varepsilon <\frac{1}{12}$ such that $0< \delta_{k} \leq \varepsilon \sigma_{k} \|\tilde{\bm{s}}_{k}\|_2^3$, then 
    \begin{equation}
        \min_{1 \leq k \leq T} \tilde{\mu}(\bm{x}_{k}) \leq \mathcal{O}\left(T^{-1/3}\right).
    \end{equation}
\end{enumerate}
\end{theorem}

\begin{proof}
Lemma \ref{lemma_5} implies that $\tilde{\mu}(\bm{x}_{k+1}) \leq \|\tilde{\bm{s}}_{k}\|_2+\frac{1}{\sqrt{c_{5}}}\left\| \nabla m_{k}(\tilde{\bm{s}}_{k}) \right\|_2^{1/2} + \frac{\sigma_{\max}}{2c_6}\left|\|\bm{s}_{k}\|_2 - \|\tilde{\bm{s}}_{k}\|_2\right|$ for all $k \in \mathcal{S}_{T}$. Then we have $\min_{1 \leq k \leq T} \tilde{\mu}(\bm{x}_{k}) \leq \min_{k \in \mathcal{S}_{T}} \|\tilde{\bm{s}}_{k}\|_2+\frac{1}{\sqrt{c_{5}}}\left\| \nabla m_{k}(\tilde{\bm{s}}_{k}) \right\|_2^{1/2} + \frac{\sigma_{\max}}{2c_6} \left|\|\bm{s}_{k}\|_2 - \|\tilde{\bm{s}}_{k}\|_2\right|$.  Combining it with Lemma \ref{lemma_7}, Lemma \ref{lemma_6} and Condition \ref{condition1}, we finish the proof. 
\end{proof}

\begin{corollary}
\label{corollary_1}
Under Assumption \ref{assump} and Condition \ref{condition1}, let the sequence $\{\bm{x}_{k}\}_{k=1}^{+\infty}$ be generated by Algorithm \ref{alg_arcm} and $f(\bm{x})$ satisfies the \KL property. If there exists $0<\varepsilon <\frac{1}{12}$ such that $0< \delta_{k} \leq \varepsilon \sigma_{k} \|\tilde{\bm{s}}_{k}\|_2^3$, then ARCm with inexact CRS solutions still hold local convergence property as in Theorem \ref{theorem_2}.
\end{corollary}

\section{Numerical Experiments}

In this section, we conduct experiments on the proposed ARCm and some state-of-the-art second-order methods (e.g., CR, CRm, ARC, and TR) in solving non-convex logistic regression and robust linear regression models.

\textbf{Settings.} All these algorithms involve solving CRS in each iteration (note that trust region subproblems are similar to cubic regularization subproblems in Step 2 of Algorithm 1). We adopt the Krylov subspace method \cite{cartis2011adaptive, conn2000trust} with at most $50$ subspaces to approximately solve CRS. All hyperparameters are tuned to achieve nearly optimal results in terms of the iterations for convergence. For ARCm, we set the starting cubic penalty parameter $\sigma_0 = 1.0$ and the momentum parameters $(\tau,\alpha_1,\alpha_2) = (0.5,0.1,1.0)$. Except for the momentum parameters, all other parameters of ARCm are the same as ARC (e.g., $(\eta_1,\eta_2)=(0.1,0.9)$).

\textbf{Datasets.} For both two models, we test algorithms on ARCENE \cite{guyon2004result}, Covtype \cite{blackard1999comparative} and DrivFace \cite{diaz2016reduced} datasets. Each training data in all three datasets is in the form of $(\bm{a}_{i},b_i)$, where $\bm{a}_i \in \mathbb{R}^{d}$ is a multi-variate attribute vector and $b_i \in \{0,1\}$ is the corresponding label. The detailed information for these two datasets is shown in Table \ref{tab1}.

\begin{table}[tbh]
\small
\center
\caption{The overview of datasets.\label{tab1}}
\begin{tabular}{lll}
 \hline
 Dataset & sample size $n$ & dimension $d$ \\
  \hline
 ARCENE & $100$ & $10000$ \\
  \hline
 Covtype & $581012$ & $54$\\
  \hline
 DrivFace & $606$ &  $6400$\\
 \hline
\end{tabular}
\end{table}

\textbf{Models.} For a given dataset $\{(\bm{a}_i,b_i)\}_{i=1}^{n}$, the empirical loss for the non-convex logistic regression is
\begin{equation*}
    \begin{split}
        \min_{\bm{w}} & \frac{1}{n} \sum_{i=1}^{n} b_i \log \left[ \psi(\bm{w}^{\top}\bm{a}_i)\right]+ (1-b_i) \log \left[ 1- \psi(\bm{w}^{\top}\bm{a}_i) \right] \\
        & + \chi \sum_{j=1}^{d}\frac{w_j^2}{1+w_j^2},
    \end{split}
\end{equation*}
where $\psi(x) = \frac{1}{1+\exp(-x)}$, $\bm{w} = (w_1~\cdots~ w_d)^{\top}$ and $\chi = 0.1$.
The empirical loss for the non-convex robust linear regression is formulated as
\begin{equation*}
    \min_{\bm{w}} \frac{1}{n} \sum_{i=1}^{n} \varrho(b_i - \bm{w}^{\top}\bm{a}_i),
\end{equation*}
where $\varrho(x) = \log(\frac{x^2}{2}+1)$.

The trajectories of gradient norm during the iterations are displayed in Figure \ref{fig_grad}, where gradient norm is the most popular measurement for convergence \cite{cartis2011adaptive, wang2020cubic}. We would like to emphasize that the Hessian is (nearly) semi-positive definite in all experiments if algorithms converge. Due to limited spaces, we put trajectories of the empirical loss in the supplement. We merely show the gradient norm versus iteration steps rather than time (seconds), since the main computation in each step is from the CRS. The proposed ARCm achieves the best performances in these examples, being $2$-$10$ times faster than CRm and $10$\%-$50$\% faster than ARC in terms of iterations for convergence, except for the logistic regression on the Covtype dataset where ARCm and ARC have similar results. We do not expect the proposed ARCm significantly outperforms ARC if the problem is less challenging and ARC can easily find the optima (Figure \ref{fig_grad} (b)). Moreover, we also highlight the two cases of momentum in Theorem \ref{theorem_0}. When the momentum term helps the convergence of ARCm (i.e., $f(\bm{z}_{k+1}) < f(\bm{y}_{k+1})$), we point it in red and gray if $\bm{s}_{k}^{\top}\bm{v}_{k-1} >0$ and $\bm{s}_{k}^{\top}\bm{v}_{k-1} <0$ respectively. We observe that the momentum in opposite direction helps the convergence of ARCm and it usually occurs at the beginning of the algorithm, which is within our expectation.

\section{Conclusion and Future Works}
In this paper, we propose the momentum accelerated ARC (ARCm) that improves the performance of ARC. The global convergence and the local convergence under the \KL property are theoretically studied. Furthermore, we analyze a more practical case where inexact CRS solutions are obtained in ARCm. Experimental results show that the proposed ARCm significantly outperforms some state-of-the-art second-order methods in solving non-convex logistic regression and robust linear regression models.  There are still many unsolved problems for cubic regularization methods. Firstly,  better strategies for CR-based and ARC-based algorithms. For example, some works modify the criterion $\rho_{k}$ in ARC. Secondly, fast solvers of CRS. Some methods enforce structured Hessian in CRS to reduce the computation.

\bibliographystyle{plain}
\bibliography{reference}

\clearpage

\appendix
\title{\Large 
Supplementary Materials}

\maketitle

\section{Some Technical Proofs}

\subsection{Proof for Lemma \ref{lemma_1}}
For any $\bm{x}, \bm{y} \in \mathcal{F}$, we have
\begin{equation}
\label{eq13}
\begin{split}
    & \left\| \nabla f(\bm{y}) - \nabla f(\bm{x}) - \nabla^2 f(\bm{x}) (\bm{y}-\bm{x}) \right\|_2\\
    & =  \left\| \int_{0}^{1} \left[\nabla^2 f(\bm{x}+t(\bm{y}-\bm{x})) - \nabla^2 f(\bm{x}) \right] (\bm{y}-\bm{x})~dt \right\|_2 \\
    & \leq L_{\rm{H}} \|\bm{y}-\bm{x}\|_2^2 \int_{0}^{1} t~dt = \frac{L_{\rm{H}}}{2} \|\bm{y}-\bm{x}\|_2^2,
\end{split}
\end{equation}
and
\begin{multline}
\label{eq15}
        \left| f(\bm{y}) - f(\bm{x}) - \nabla f(\bm{x}_{k})^{\top} (\bm{y}-\bm{x}) \right.\\
         - \left. \frac{1}{2} (\bm{y}-\bm{x})^{\top} \nabla^2 f(\bm{x}) (\bm{y}-\bm{x}) \right| \\
        = \left| \int_{0}^{1} \left( \nabla f(\bm{x}+t (\bm{y}-\bm{x})) - \nabla f(\bm{x}) \right.\right.\\
        \left.\left. - t \nabla^2 f(\bm{x}) (\bm{y}-\bm{x}) \right)^{\top}(\bm{y}-\bm{x})~ dt \right|\\
        \leq \frac{L_{\rm{H}}}{2} \|\bm{y}-\bm{x}\|_2^3 \cdot \int_{0}^{1} t^2~dt = \frac{L_{\rm{H}}}{6} \|\bm{y}-\bm{x}\|_2^3.  
\end{multline}

Then
\begin{equation*}
    f(\bm{x}_{k} + \bm{s}_{k}) - m_k(\bm{s}_{k}) \leq \left( \frac{L_{\rm{H}}}{6} - \frac{\sigma_k}{6} \right) \|\bm{s}_{k}\|_2^3.
\end{equation*}
Let $r_k := f(\bm{x}_{k} + \bm{s}_{k}) - m_k(\bm{s}_{k}) + (1-\eta_2) \left(m_{k}(\bm{s}_{k}) - f(\bm{x}_{k}) \right)$. The very successful update occurs (i.e., $\rho_k > \eta_2$) if and only if $r_k<0$. Note that $m_{k}(\bm{s}_{k}) - f(\bm{x}_{k}) < 0$ if $\bm{x}_{k}$ is not a local minimizer. If $\sigma_{k} \geq L_{\rm{H}}$, then $f(\bm{x}_{k} + \bm{s}_{k}) - m_k(\bm{s}_{k}) \leq \left( \frac{L_{\rm{H}}}{6} - \frac{\sigma_k}{6} \right) \|\bm{s}_{k}\|_2^3 \leq 0$ and $r_{k} < 0$. Suppose that $k-1 \notin \mathcal{S}$ and $\sigma_{k-1} \approx L_{\rm{H}}$ (but $\sigma_{k-1} < L_{\rm{H}}$), then $L_{\rm{H}} < \sigma_{k} = \gamma_1 \sigma_{k-1} \leq \gamma_1 L_{\rm{H}}$ and $k \in \mathcal{S}$.

\subsection{Proof for Lemma \ref{lemma_8}}
Note that the sequence $\{f(\bm{x}_{k})\}$ is lower bounded by $f^{*}$ and is monotonically decreasing, then it must be convergent. Assumption \ref{assump2} and the non-increasing property of $\{f(\bm{x}_{k})\}$ imply the boundedness of $\{\bm{x}_{k}\}_{k=0}^{+\infty}$. Here the boundedness of the sequence $\{\bm{x}_{k}\}_{k=0}^{+\infty}$ implies that the set $\bar{\mathcal{X}}$ is non-empty. Using Lemma \ref{lemma_4}, we have for all $k \in \mathcal{S}$
\begin{equation*}
    \begin{split}
        & \|\bm{x}_{k+1} - \bm{x}_{k}\|_2\\
        \leq & \|\bm{x}_{k+1} - \bm{y}_{k+1}\|_2 + \|\bm{y}_{k+1} - \bm{x}_{k}\|_2\\
        \leq & \beta_{k} \|\bm{v}_{k}\|_2 + \|\bm{s}_{k}\|_2\\
        \leq & \alpha_1 \|\bm{s}_{k}\|_2 \cdot  \|\bm{v}_{k}\|_2 + \|\bm{s}_{k}\|_2\\
        \leq & \left( \alpha_1 \frac{1}{1-\tau}  \left( \frac{12(f(\bm{x}_0 ) - f^{*})}{\eta_1 \sigma_{\min}} \right)^{1/3} + 1 \right) \|\bm{s}_{k}\|_2\\
        := & c_{7} \|\bm{s}_{k}\|_2.
    \end{split}
\end{equation*}
Lemma \ref{lemma_3} and \ref{lemma_4} jointly show that $\lim_{k \to \infty, k \in \mathcal{S}} \|\bm{s}_{k}\|_2 = 0$ and thus $\lim_{k \to \infty}\|\bm{x}_{k+1} - \bm{x}_{k}\|_2 = 0$ (note that $\|\bm{x}_{k+1} - \bm{x}_{k}\|_2 = 0$ if $k \notin \mathcal{S}$). Using Lemma \ref{lemma_2}, we have
\begin{equation*}
\begin{split}
    \left\|\nabla f(\bar{\bm{x}})\right\| & \leq \limsup_{k \to \infty, k \in \mathcal{S}} \left\|\nabla f(\bm{x}_{k+1})\right\| \\
    & \leq \limsup_{k \to \infty, k \in \mathcal{S}} c_1 \|\bm{s}_{k}\|_2^2 = 0
\end{split}
\end{equation*}
and
\begin{equation*}
\begin{split}
    \lambda_{\min}\left(\nabla^2 f(\bar{\bm{x}})\right) & \geq \limsup_{k \to \infty, k \in \mathcal{S}}  \lambda_{\min}\left(\nabla^2 f(\bm{x}_{k+1})\right) \\
    & \geq \limsup_{k \to \infty, k \in \mathcal{S}}  -c_2\|\bm{s}_{k}\|_2 = 0.
\end{split}
\end{equation*}

Define $e_{j} = f(\bm{x}_{k_j+1}) - \bar{f}$, where $k_j \in \mathcal{S}$ is the $j$-th element in $\mathcal{S}$. If $f(\cdot)$ satisfies the \KL property, then there exist a large enough $j_0$ such that $\phi^{\prime}(e_{j}) \geq \frac{1}{\left\| \nabla f(\bm{x}_{k_j+1}) \right\|}$ and $1-\alpha_2 \|\bm{s}_{k_{j}}\|_2 \|\bm{v}_{k_j}\|_2 \geq \alpha_3$  for all $j \geq j_0$ and some $1>\alpha_3>0$, since $\lim_{j \to \infty} \|\bm{s}_{k_j}\|_2=0$ and $\|\bm{v}_{k_j}\|_2$ is bounded. Furthermore,
\begin{equation}
\label{eq6}
\begin{split}
    \|\bm{x}_{k_j+1} - \bm{x}_{k_j}\|_2 & \geq \|\bm{s}_{k_j}\|_2 - \beta_{k_j}\|\bm{v}_{k_j}\|_2 \\
    & \geq \|\bm{s}_{k_j}\|_2 - \alpha_2 \|\bm{s}_{k_j}\|_2^2 \|\bm{v}_{k_j}\|_2\\
    & = \left(1-\alpha_2 \|\bm{s}_{k_j}\|_2 \|\bm{v}_{k_j}\|_2 \right) \|\bm{s}_{k_j}\|_2\\
    & \geq \alpha_3 \|\bm{s}_{k_j}\|_2,
\end{split}
\end{equation}

\begin{equation*}
\begin{split}
    \|\bm{x}_{k_j+1} - \bm{x}_{k_j}\|_2 & \leq \|\bm{s}_{k_j}\|_2 + \beta_{k_j}\|\bm{v}_{k_j}\|_2\\
    & \leq \|\bm{s}_{k_j}\|_2 + \alpha_2 \|\bm{s}_{k_j}\|_2^2 \|\bm{v}_{k_j}\|_2\\
    &  = \left(1+\alpha_2 \|\bm{s}_{k_j}\|_2 \|\bm{v}_{k_j}\|_2 \right) \|\bm{s}_{k_j}\|_2\\
    & \leq (2-\alpha_3) \|\bm{s}_{k_j}\|_2,
\end{split}
\end{equation*}

and
\begin{equation*}
    \phi^{\prime}(e_{j}) \geq \frac{1}{\left\| \nabla f(\bm{x}_{k_j+1}) \right\|} \geq \frac{1}{c_1 \|\bm{s}_{k_j}\|_2^2} \geq \frac{\alpha_3^2}{c_1\|\bm{x}_{k_j+1} - \bm{x}_{k_j}\|_2^2}.
\end{equation*}
Note that $\phi(\cdot)$ is concave when $\theta < 1$ and $e_{j+1} < e_{j}$, we have
\begin{equation}
\label{eq10}
\begin{split}
    & \phi(e_{j}) - \phi(e_{j+1})\\
    \geq & \phi^{\prime}(e_{j}) (e_{j} - e_{j+1})\\
    \geq & \frac{\alpha_3^2}{c_1\|\bm{x}_{k_j+1} - \bm{x}_{k_j}\|_2^2} \cdot \frac{1}{12}\eta_1 \sigma_{\min}\|\bm{s}_{k_{j+1}}\|_2^3\\
    \geq & \frac{\alpha_3^2 \eta_1 \sigma_{\min}}{12c_1 (2-\alpha_3)} \frac{\|\bm{x}_{k_{j+1}+1} - \bm{x}_{k_{j+1}}\|_2^3}{\|\bm{x}_{k_j+1} - \bm{x}_{k_j}\|_2^2}\\
    := & c_{8} \frac{\|\bm{x}_{k_{j+1}+1} - \bm{x}_{k_{j+1}}\|_2^3}{\|\bm{x}_{k_j+1} - \bm{x}_{k_j}\|_2^2}.
\end{split}
\end{equation}
Then by the H\"{o}lder's inequality, we have 
\begin{equation*}
\label{eq5}
\begin{split}
    & \sum_{j=j_0}^{J}  \|\bm{x}_{k_{j+1}+1} - \bm{x}_{k_{j+1}}\|_2\\
    & \leq c_{9} \cdot \sum_{j=j_0}^{J} (\phi(e_{j}) - \phi(e_{j+1}))^{1/3} \|\bm{x}_{k_j+1} - \bm{x}_{k_j}\|_2^{2/3}\\
    & \leq c_{9} \left(\sum_{j=j_0}^{J}\phi(e_{j}) - \phi(e_{j+1}) \right)^{1/3} \left(\sum_{j=j_0}^{J} \|\bm{x}_{k_j+1} - \bm{x}_{k_j}\|_2 \right)^{2/3}\\
    & \leq c_{9} \phi(e_{j_0})^{1/3} \cdot \left(\sum_{j=j_0}^{J} \|\bm{x}_{k_j+1} - \bm{x}_{k_j}\|_2 \right)^{2/3},
\end{split}
\end{equation*}
where $c_{9} = c_{8}^{-1/3}$. 
If $\sum_{j=j_0}^{+\infty}  \|\bm{x}_{k_{j}+1} - \bm{x}_{k_{j}}\|_2 = + \infty$, then we may let $J>j_0$ to be large enough such that $\sum_{j=j_0}^{J} \|\bm{x}_{k_j+1} - \bm{x}_{k_j}\|_2  \approx \sum_{j=j_0}^{J}  \|\bm{x}_{k_{j+1}+1} - \bm{x}_{k_{j+1}}\|_2 \gg 1$. However the relation (\ref{eq5}) is violated. Therefore  $\sum_{j=j_0}^{+\infty}  \|\bm{x}_{k_{j}+1} - \bm{x}_{k_{j}}\|_2 < + \infty$ and thus $\sum_{k=0}^{+\infty}  \|\bm{x}_{k+1} - \bm{x}_{k}\|_2 < + \infty$, which implies that $\{\bm{x}_{k}\}_{k=0}^{+\infty}$ is a Cauchy sequence and $\bar{\mathcal{X}}=\{\bar{\bm{x}}\}$ is a singleton.

\subsection{Proof for Theorem \ref{theorem_2}}

First note that if $f(\bm{x})$ satisfies the \KL property, then it also satisfies the local \KL-error bound \cite[Proposition~1]{zhou2018convergence}, which is a generalization of the local error bound in \cite{yue2019quadratic}, i.e., 
\begin{equation}
\label{eq7}
    \|\bm{x}_{k_{j}+1} - \bar{\bm{x}}\|_2 \leq \kappa \left\|\nabla f(\bm{x}_{k_{j}+1})\right\|_2^{\frac{\theta}{1-\theta}}
\end{equation}
holds for $j>j_0$ with large enough $j_0>0$ and some $\kappa > 0$. Combining (\ref{eq7}) with (\ref{eq6}) and Lemma \ref{lemma_2}, we have
\begin{equation}
\label{eq8}
    \begin{split}
        \|\bm{x}_{k_{j}+1} - \bar{\bm{x}}\|_2 & \leq \kappa \left\|\nabla f(\bm{x}_{k_{j}+1})\right\|_2^{\frac{\theta}{1-\theta}}\\
        & \leq \kappa \left(c_1 \|\bm{s}_{k_j}\|_2^2 \right)^{\frac{\theta}{1-\theta}}
    \end{split}
\end{equation}
Note that \cite[Lemma~1]{yue2019quadratic} shows that there exists $c_{10}>0$ such that 
\begin{equation}
\label{eq9}
    \|\bm{s}_{k_j}\|_2 \leq c_{10} \|\bm{x}_{k_{j}} - \bar{\bm{x}}\|_2,
\end{equation}
where $c_{10} = \left(1+\frac{L_{\rm{H}}}{\sigma_{\min}} + \sqrt{\left(1+\frac{L_{\rm{H}}}{\sigma_{\min}} \right)^2 + \frac{L_{\rm{H}}}{\sigma_{\min}}} \right)$. Equations (\ref{eq8}) and (\ref{eq9}) jointly imply that 
\begin{equation*}
    \|\bm{x}_{k_{j}+1} - \bar{\bm{x}}\|_2 \leq c_{11} \left\|\bm{x}_{k_{j}} - \bar{\bm{x}}\right\|_2^{\frac{2\theta}{1-\theta}} = c_{11} \left\|\bm{x}_{k_{j-1}+1} - \bar{\bm{x}}\right\|_2^{\frac{2\theta}{1-\theta}},
\end{equation*}
where the convergence is super linear since $\frac{2\theta}{1-\theta}>1$ when $\theta \in (\frac{1}{3},1)$. Without the loss of generality we assume that $c_{11} \leq 1$ and $\|\bm{x}_{k_{j_0}+1} - \bar{\bm{x}}\|_2 \leq \exp(-1)$, we conclude the results in item 1 for $\theta \in (\frac{1}{3},1)$.

We then discuss the cases when $\theta=\frac{1}{3}$ and $\theta \in (0,\frac{1}{3})$. Here we adopt similar techniques in \cite[Theorem~4]{zhou2018convergence}. Fix $\nu \in (0,1)$ and $j>j_0$ for a large enough $j_0 \in \mathcal{S}$. If $\|\bm{x}_{k_{j+1}+1} - \bm{x}_{k_{j+1}}\|_2 \geq \nu \|\bm{x}_{k_{j}+1} - \bm{x}_{k_{j}}\|_2$, then using (\ref{eq10}) we have
\begin{equation*}
\begin{split}
    \phi(e_{j}) - \phi(e_{j+1})  & \geq c_{8} \frac{\|\bm{x}_{k_{j+1}+1} - \bm{x}_{k_{j+1}}\|_2^3}{\|\bm{x}_{k_j+1} - \bm{x}_{k_j}\|_2^2}\\
    & \geq c_{8} \nu^2 \|\bm{x}_{k_{j+1}+1} - \bm{x}_{k_{j+1}}\|_2
\end{split}
\end{equation*}
and equivalently
\begin{equation*}
    \|\bm{x}_{k_{j+1}+1} - \bm{x}_{k_{j+1}}\|_2 \leq \frac{1}{c_{8} \nu^2} \left( \phi(e_{j}) - \phi(e_{j+1}) \right).
\end{equation*}
Otherwise, $\|\bm{x}_{k_{j+1}+1} - \bm{x}_{k_{j+1}}\|_2 \leq \nu \|\bm{x}_{k_{j}+1} - \bm{x}_{k_{j}}\|_2$. Summing up two inequalities, we have
\begin{equation}
\label{eq11}
    \begin{split}
        & \|\bm{x}_{k_{j+1}+1} - \bm{x}_{k_{j+1}}\|_2\\
        & \leq \frac{1}{c_{8} \nu^2} \left( \phi(e_{j}) - \phi(e_{j+1}) \right) + \nu \|\bm{x}_{k_{j}+1} - \bm{x}_{k_{j}}\|_2.
    \end{split}
\end{equation}
Let $\Delta_{j} = \sum_{l=j}^{\infty} \|\bm{x}_{k_{l}+1} - \bm{x}_{k_{l}}\|_2$, then $\|\bm{x}_{k_{j}+1} - \bm{x}_{k_{j}}\|_2 = \Delta_{j}-\Delta_{j+1}$. Summing up (\ref{eq11}) from $j_0$ to $J$ yields that
\begin{equation*}
\begin{split}
    & \sum_{j=j_0}^{J}\|\bm{x}_{k_{j+1}+1} - \bm{x}_{k_{j+1}}\|_2\\
    & \leq \nu \sum_{j=j_0}^{J}\|\bm{x}_{k_{j}+1} - \bm{x}_{k_{j}}\|_2 + \frac{1}{c_{8} \nu^2} \phi(e_{j_0}).
\end{split}
\end{equation*}
and then 
\begin{equation*}
\begin{split}
    & \sum_{j=j_0}^{J}\|\bm{x}_{k_{j+1}+1} - \bm{x}_{k_{j+1}}\|_2 \\
    \leq & \nu \left(\sum_{j=j_0}^{J}\|\bm{x}_{k_{j+1}+1} - \bm{x}_{k_{j+1}}\|_2 + \|\bm{x}_{k_{j_0}+1} - \bm{x}_{k_{j_0}}\|_2 \right)\\
    & + \frac{1}{c_{8} \nu^2} \phi(e_{j_0}).
\end{split}
\end{equation*}

Let $J \to \infty$ and we can further simplify the above inequality as
\begin{equation*}
\small
\begin{split}
    & \Delta_{j+1}\\
    \leq & \frac{\nu}{1-\nu} \|\bm{x}_{k_{j}+1} - \bm{x}_{k_{j}}\|_2 + \frac{1}{c_{8}\nu^2(1-\nu)} \phi(e_{j})\\
    = & \frac{\nu}{1-\nu} \|\bm{x}_{k_{j}+1} - \bm{x}_{k_{j}}\|_2 + \frac{c}{c_{8} \nu^2(1-\nu)\theta} e_j^{\theta}\\
    \leq & \frac{\nu}{1-\nu} \|\bm{x}_{k_{j}+1} - \bm{x}_{k_{j}}\|_2 + \frac{c}{c_{8} \nu^2(1-\nu)\theta} \left(c_0 \left\| \nabla f(\bm{x}_{k_j+1}) \right\|^{\frac{1}{1-\theta}} \right)^{\theta}\\
    \leq & \frac{\nu}{1-\nu} \|\bm{x}_{k_{j}+1} - \bm{x}_{k_{j}}\|_2 + \frac{c}{c_{8} \nu^2(1-\nu)\theta} \left(c_0 \left(c_1\left\|\bm{s}_{k_j} \right\|^2\right)^{\frac{1}{1-\theta}} \right)^{\theta}\\
    \leq & \frac{\nu}{1-\nu} \|\bm{x}_{k_{j}+1} - \bm{x}_{k_{j}}\|_2 +  c_{12} \left\|\bm{x}_{k_{j}+1} - \bm{x}_{k_{j}}\right\|_2^{\frac{2\theta}{1-\theta}}\\
    = & \frac{\nu}{1-\nu} \left( \Delta_{j}-\Delta_{j+1}\right) + c_{12} \left( \Delta_{j}-\Delta_{j+1}\right)^{\frac{2\theta}{1-\theta}}.
\end{split}
\end{equation*}
If $\theta = \frac{1}{3}$, then $\Delta_{j+1} \leq c_{13}\left( \Delta_{j}-\Delta_{j+1}\right)$ and thus $\Delta_{j+1}  \leq \frac{c_{13}}{1+c_{13}} \Delta_{j}$. Therefore, $\Delta_{j} \leq \left( \frac{c_{13}}{1+c_{13}}\right)^{j-j_0} \Delta_{j_0}$ and $\|\bm{x}_{k_j+1} - \bar{\bm{x}}\|_2 \leq \Delta_{j+1} = \mathcal{O}\left(\exp \left(-c_{14}(j-j_0) \right) \right)$, where $c_{14} = \log\left(\frac{1+c_{13}}{c_{13}}\right)$. If $\theta \in (0,\frac{1}{3})$, then $\left( \Delta_{j}-\Delta_{j+1}\right)^{\frac{2\theta}{1-\theta}}$ is dominant term of the right hand side if $j_0$ is large enough. Then there exists $c_{15}>0$ such that $\Delta_{j+1}^{\frac{1-\theta}{2\theta}} \leq c_{15}\left( \Delta_{j}-\Delta_{j+1}\right)$. Define $h(t) = t^{-\frac{1-\theta}{2\theta}}$ and a fixed constant $\omega > 1$. If $h(\Delta_{k+1}) \leq \omega h(\Delta_{k})$, then 
\begin{equation*}
    \begin{split}
        1 & \leq c_{15}\left( \Delta_{j}-\Delta_{j+1}\right) h(\Delta_{j+1}) \leq c_{15} \omega \left( \Delta_{j}-\Delta_{j+1}\right)h(\Delta_{j}) \\
        & \leq c_{15} \omega \int_{\Delta_{j+1}}^{\Delta_{j}} h(t)~dt = c_{15} \omega \frac{2\theta}{3\theta -1} \left(\Delta_{j}^{\frac{3\theta-1}{2\theta}} - \Delta_{j+1}^{\frac{3\theta-1}{2\theta}} \right)
    \end{split}
\end{equation*}
and
\begin{equation}
\label{eq12}
     \Delta_{j+1}^{\frac{3\theta-1}{2\theta}} - \Delta_{j}^{\frac{3\theta-1}{2\theta}} \geq \frac{1-3\theta}{2\omega \theta c_{15}}.
\end{equation}
On the other hand, if $h(\Delta_{k+1}) > \omega h(\Delta_{k})$, then $\Delta_{j+1}^{\frac{3\theta-1}{2\theta}} > \omega^{\frac{1-3\theta}{1-\theta}}\Delta_{j}^{\frac{3\theta-1}{2\theta}}$ and $\Delta_{j+1}^{\frac{3\theta-1}{2\theta}} - \Delta_{j}^{\frac{3\theta-1}{2\theta}}  > \left(\omega^{\frac{1-3\theta}{1-\theta}}-1\right)\Delta_{j}^{\frac{3\theta-1}{2\theta}}$. Note that $\omega^{\frac{1-3\theta}{1-\theta}}-1 > 0$ and $\lim_{j \to \infty} \Delta_{j}^{\frac{3\theta-1}{2\theta}} = + \infty$, then there must exist large enough $j_0$ such that $\left(\omega^{\frac{1-3\theta}{1-\theta}}-1\right)\Delta_{j}^{\frac{3\theta-1}{2\theta}} \geq \frac{1-3\theta}{2\omega \theta c_{15}}$. Therefore, the inequality (\ref{eq12}) holds for all $j>j_0$. Summing up (\ref{eq12}) from $j=j_0$ we have that 
\begin{equation*}
    \Delta_{j+1}^{\frac{3\theta-1}{2\theta}} - \Delta_{j_{0}}^{\frac{3\theta-1}{2\theta}} \geq \frac{1-3\theta}{2\omega \theta c_{15}}(j+1-j_0)
\end{equation*}
and 
\begin{equation*}
\begin{split}
    \Delta_{j+1} & \leq \left(\Delta_{j_{0}}^{\frac{3\theta-1}{2\theta}} + \frac{1-3\theta}{2\omega \theta c_{15}}(j+1-j_0) \right)^{-\frac{2\theta}{1-3\theta}}\\
    & = \mathcal{O}\left(\left(c_{16} \left(j-j_0\right) \right)^{-\frac{2\theta}{1-3\theta}}\right),
\end{split}
\end{equation*}
which completes the proof since $\|\bm{x}_{k_j+1} - \bar{\bm{x}}\|_2 \leq \Delta_{j+1}$.

\subsection{Proof for Lemma \ref{lemma_7}}
A similar analysis is conducted here as Lemma \ref{lemma_1}. Let $\tilde{r}_{k} = f(\bm{x}_{k}+\tilde{\bm{s}}_{k}) - m_{k}(\tilde{\bm{s}}_{k})+ (1-\eta_2)(m_{k}(\tilde{\bm{s}}_{k})-f(\bm{x}_{k}))$. According to (\ref{eq13}), we have
\begin{equation*}
    f(\bm{x}_{k}+\tilde{\bm{s}}_{k}) - m_{k}(\tilde{\bm{s}}_{k}) \leq \left(\frac{L_{\rm{H}}}{6}-\frac{\sigma_{k}}{6} \right)\left\|\tilde{\bm{s}}_{k}\right\|_2^3.
\end{equation*}
The proof can be similarly developed by using the argument in Lemma \ref{lemma_1}.

\subsection{Proof for Lemma \ref{lemma_6}}
As an analogy to Lemma \ref{lemma_4}, we provide only the key steps here and omit some details. The inequality (\ref{eq14}) can be similarly developed as
\begin{equation*}
    \sum_{k \in \mathcal{S}_{T}} \eta_1\left(\frac{\sigma_{k}}{12}\left\|\tilde{\bm{s}}_{k}\right\|_2^3 - \delta_{k} \right) \leq f(\bm{x}_{0}) - f^{*},
\end{equation*}
then 
\begin{equation*}
    \max_{k \in \mathcal{S}_{T}} \frac{\sigma_{k}}{12}\left\|\tilde{\bm{s}}_{k}\right\|_2^3 \leq \frac{f(\bm{x}_{0}) - f^{*}}{\eta_1} + \delta_{k} \leq \frac{2\left(f(\bm{x}_{0}) - f^{*} \right)}{\eta_1}
\end{equation*}
and
\begin{equation*}
    \min_{k \in \mathcal{S}_{T}} \frac{\sigma_{k}}{12}\left\|\tilde{\bm{s}}_{k}\right\|_2^3 -  \delta_{k} \leq \frac{f(\bm{x}_{0}) - f^{*} }{\eta_1\left|\mathcal{S}_{T} \right|}.
\end{equation*}

\begin{figure*}[t!]
  \centering
    \subfloat[ARCENE]{
    \label{logistic_loss_arcene}
    \includegraphics[scale=0.3]{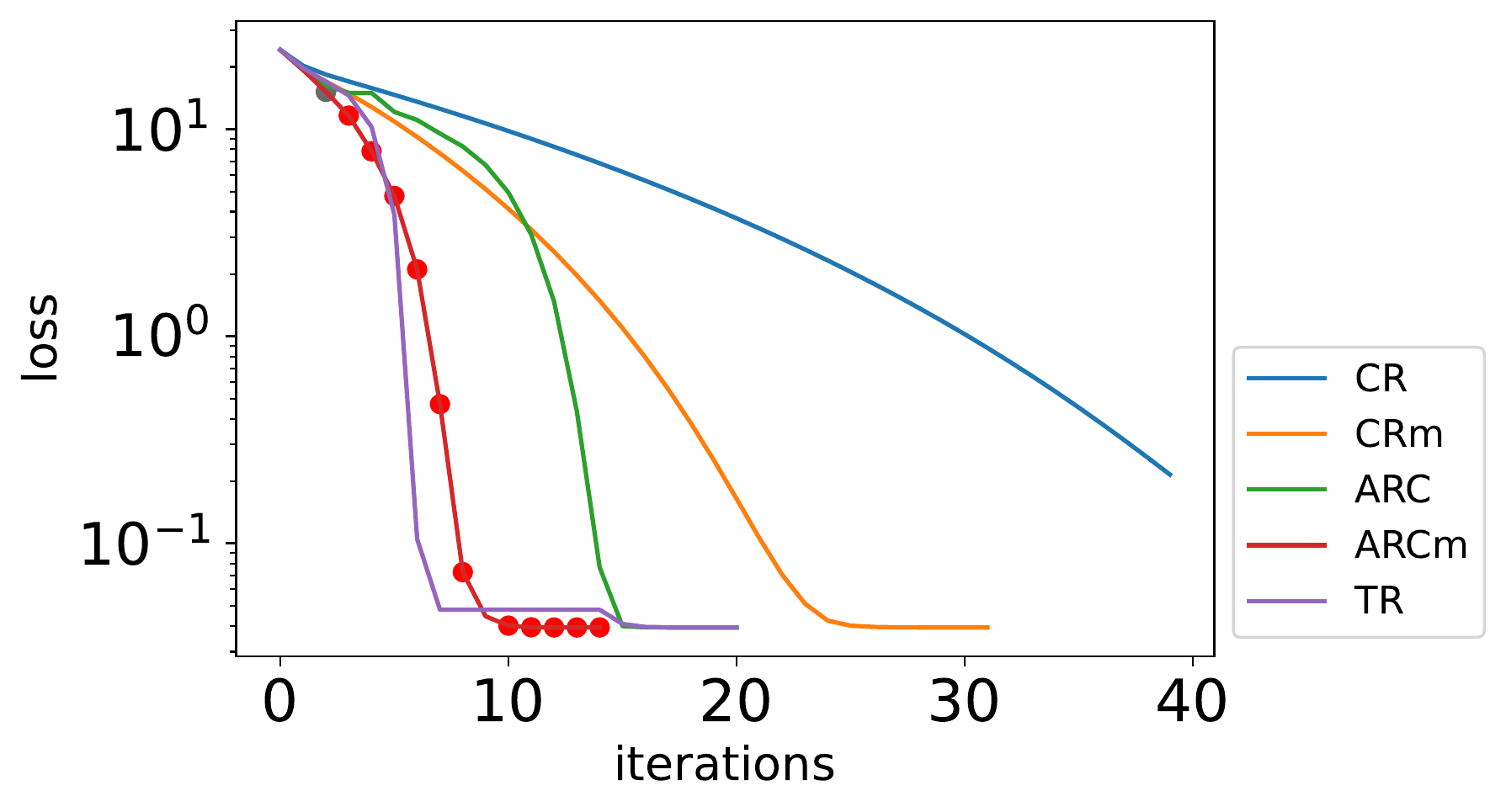}
    }\
    \subfloat[Covtype]{
    \label{logistic_loss_covtype}
    \includegraphics[scale=0.3]{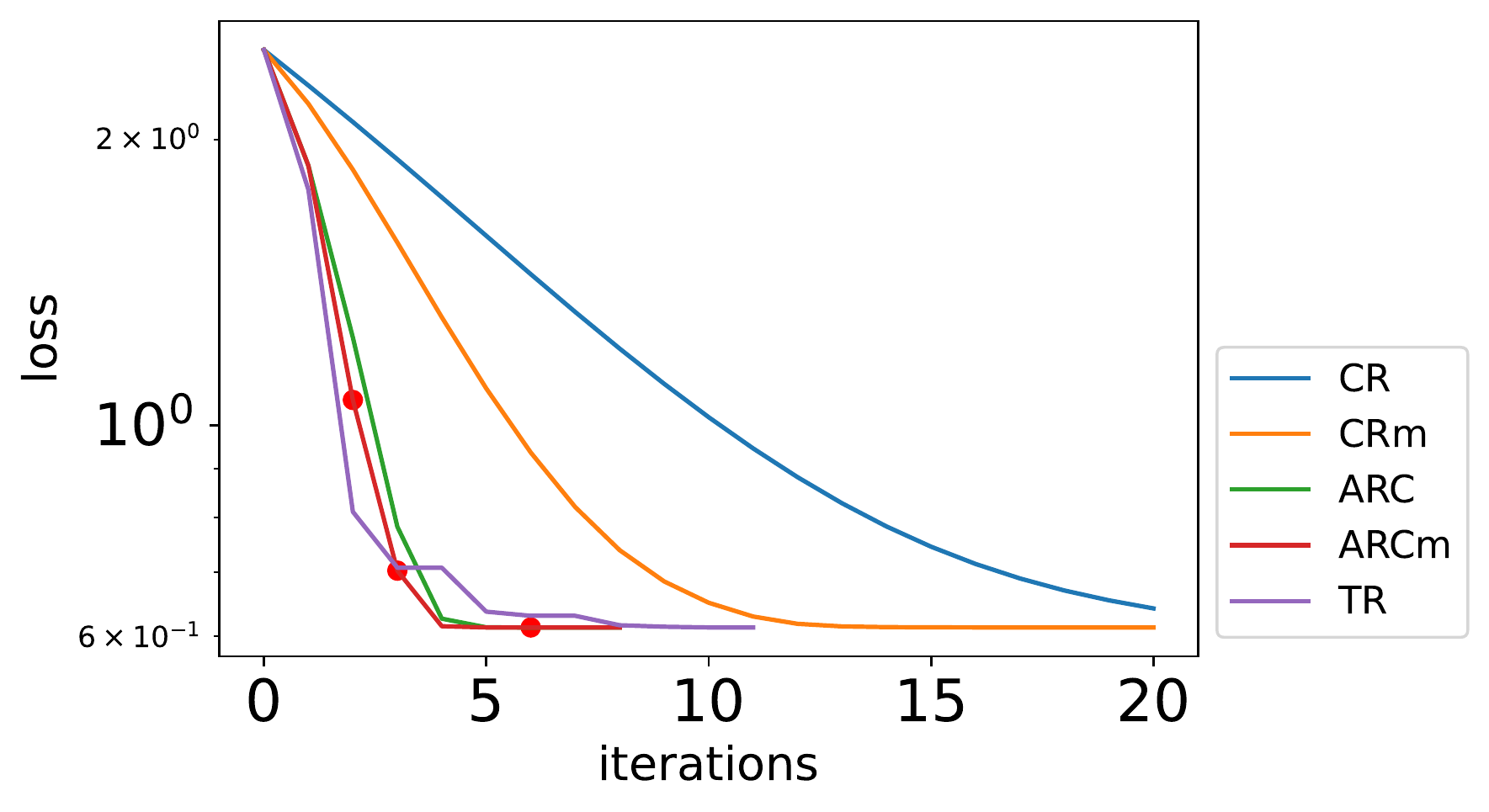}
    }\
    \subfloat[DrivFace]{
    \label{logistic_loss_drivface}
    \includegraphics[scale=0.3]{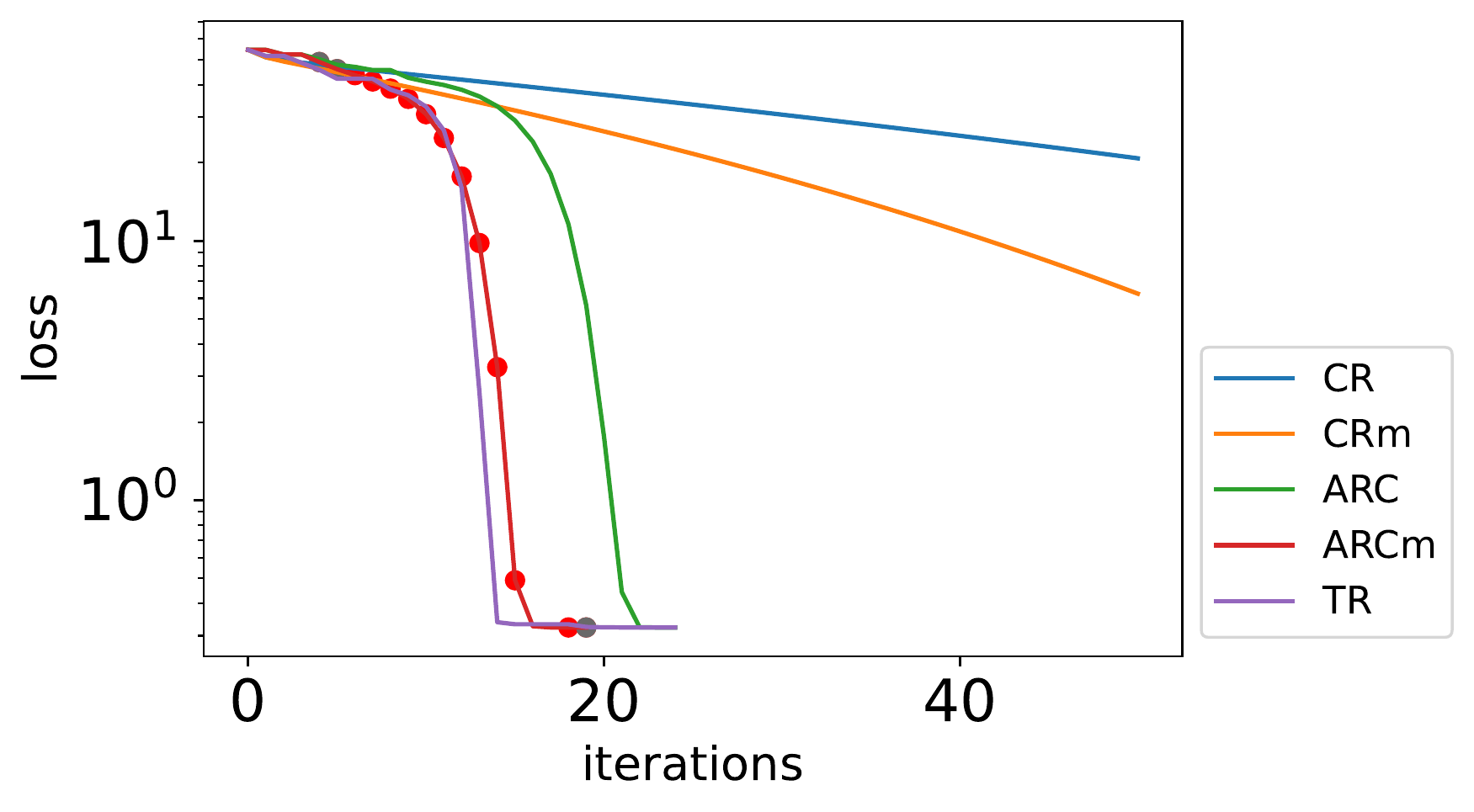}
    }\\
    \subfloat[ARCENE]{
    \label{robust_loss_arcene}
    \includegraphics[scale=0.3]{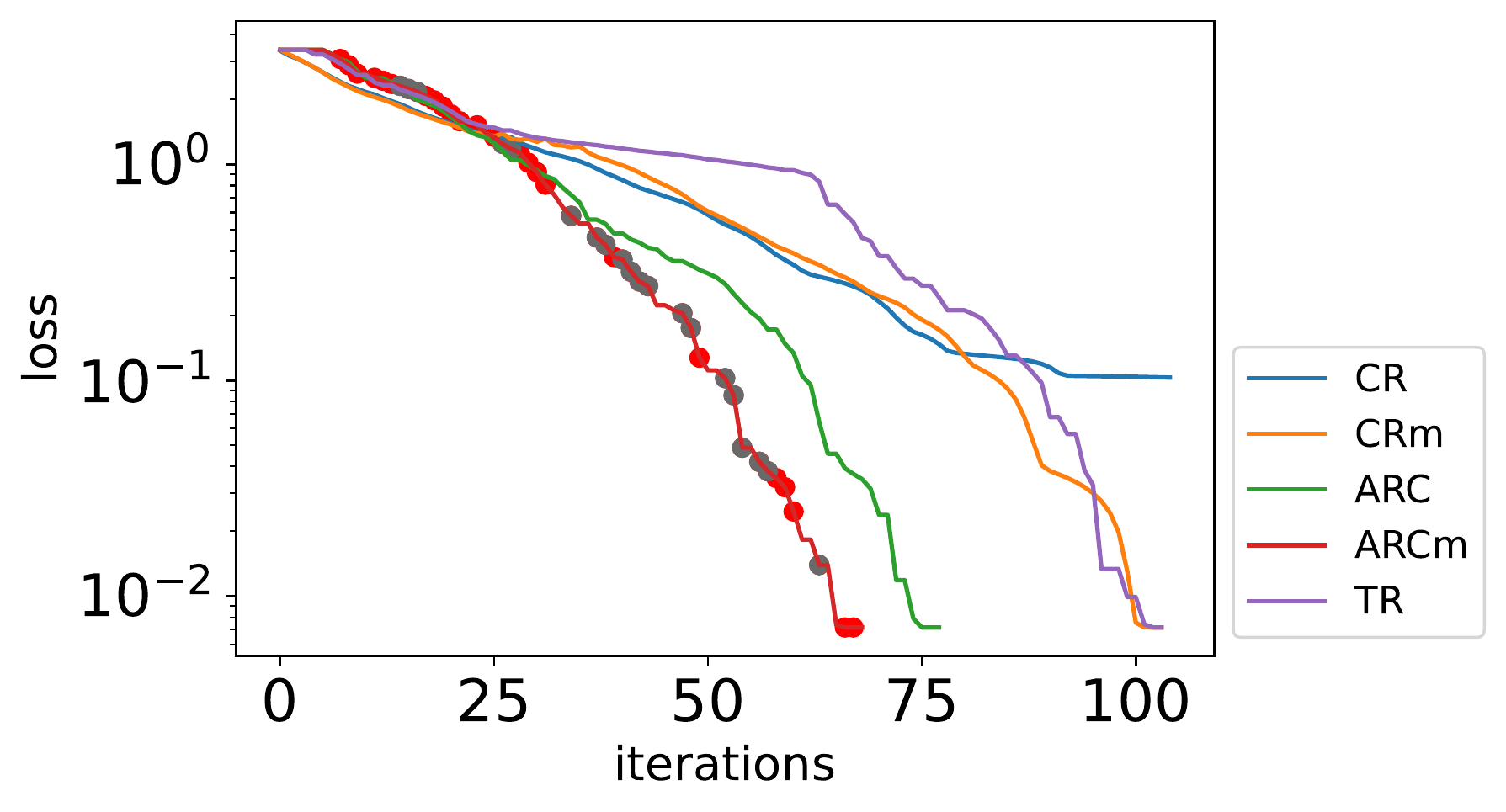}
    }\
    \subfloat[Covtype]{
    \label{robust_loss_covtype}
    \includegraphics[scale=0.3]{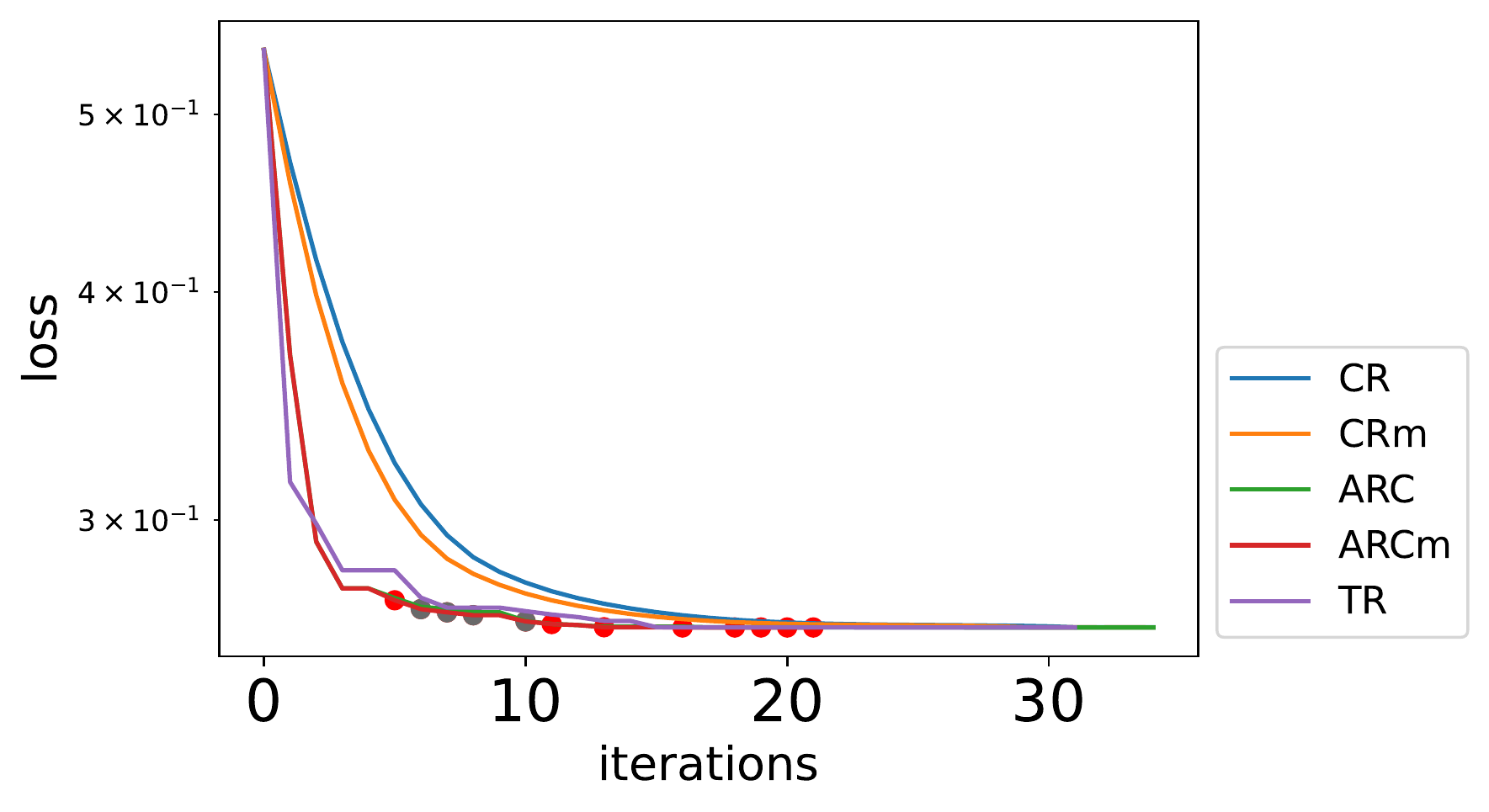}
    }\
    \subfloat[DrivFace]{
    \label{robust_loss_drivface}
    \includegraphics[scale=0.3]{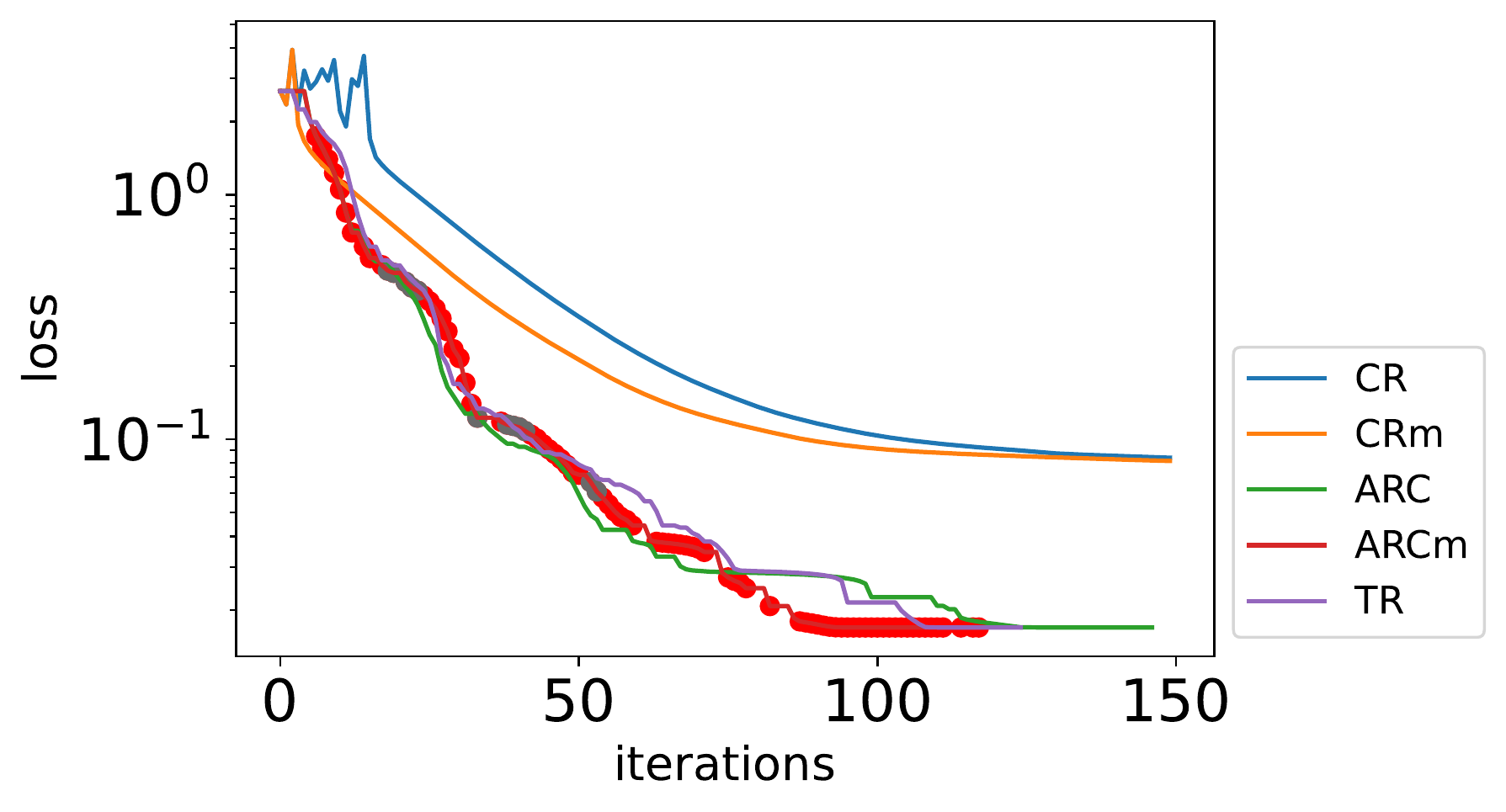}
    }\
    \caption{Loss versus the iteration steps on logistic regression (the first row) and robust linear regression (the second row) models for three datasets. Red dots: $\beta_{k}>0$ and $\bm{s}_{k}^{\top}\bm{v}_{k-1} >0$; gray dots: $\beta_{k}>0$ and $\bm{s}_{k}^{\top}\bm{v}_{k-1} < 0$.}
    \label{fig_loss}
\end{figure*}

\subsection{Proof for Lemma \ref{lemma_5}}

We first derive the error bound for $\left\|\nabla f(\bm{y}_{k+1}) \right\|$ and $\lambda_{\min}\left(\nabla^2 f(\bm{y}_{k+1}) \right)$: 

\begin{equation}
    \begin{split}
        & \left\| \nabla f(\bm{y}_{k+1}) \right\|_2\\
        \leq & \left\| \nabla f(\bm{y}_{k+1}) - \nabla f(\bm{x}_{k}) - \nabla^2 f(\bm{x}_{k}) \tilde{\bm{s}}_{k} \right\|_2\\
        & + \left\| \nabla f(\bm{x}_{k}) + \nabla^2 f(\bm{x}_{k}) \tilde{\bm{s}}_{k} + \frac{\sigma_{k}}{2}\|\tilde{\bm{s}}_{k}\|_2 \tilde{\bm{s}}_{k} \right\|_2 + \frac{\sigma_{k}}{2}\|\tilde{\bm{s}}_{k}\|_2^2\\
        \leq & \frac{L_{\rm{H}}}{2}\|\tilde{\bm{s}}_{k}\|_2^2 + \left\| \nabla m_{k}(\tilde{\bm{s}}_{k}) \right\|_2 + \frac{\sigma_{k}}{2}\|\tilde{\bm{s}}_{k}\|_2^2\\
        \leq & \left(\frac{\sigma_{\max}}{2} +\frac{L_{\rm{H}}}{2}   \right)\|\tilde{\bm{s}}_{k}\|_2^2 + \left\| \nabla m_{k}(\tilde{\bm{s}}_{k}) \right\|_2,
    \end{split}
\end{equation}
where the second inequality is due to (\ref{eq13}); and

\begin{equation*}
\begin{split}
    & \lambda_{\min}\left(\nabla^2 f(\bm{y}_{k+1}) \right)\\
    \geq & \lambda_{\min}\left(\nabla^2 f(\bm{x}_{k}) \right) - \left\|\nabla^2 f(\bm{y}_{k+1}) - \nabla^2 f(\bm{x}_{k}) \right\|_2\\
    \geq & -\frac{\sigma_{k}}{2}\|\bm{s}_{k}\|_2 - L_{\rm{H}} \|\tilde{\bm{s}}_{k}\|_2\\
    \geq & -\frac{\sigma_k}{2} \left|\|\bm{s}_{k}\|_2 - \|\tilde{\bm{s}}_{k}\|_2\right| - \frac{\sigma_k}{2} \|\tilde{\bm{s}}_{k}\|_2 - L_{\rm{H}} \|\tilde{\bm{s}}_{k}\|_2\\
    = & -\left(\frac{\sigma_{\max}}{2} - L_{\rm{H}} \right) \left\|\tilde{\bm{s}}_{k}\right\|_2 - \frac{\sigma_{\max}}{2} \left|\|\bm{s}_{k}\|_2 - \|\tilde{\bm{s}}_{k}\|_2\right|.
\end{split}
\end{equation*}
where the second inequality comes from a well-known result of cubic regularization \cite[Proposition~1]{nesterov2006cubic}. Now, we are ready to develop the error bound $\left\|\nabla f(\bm{x}_{k+1}) \right\|$ and $\lambda_{\min}\left(\nabla^2 f(\bm{x}_{k+1}) \right)$. Using similar arguments as in Lemma \ref{lemma_2}, we have 
\begin{equation*}
    \begin{split}
        & \left \| \nabla f(\bm{x}_{k+1}) \right \|_2\\
        \leq & \left \| \nabla f(\bm{y}_{k+1}) \right \|_2 + \left \| \nabla f(\bm{x}_{k+1}) - \nabla f(\bm{y}_{k+1}) \right \|_2\\
        \leq & \left \| \nabla f(\bm{y}_{k+1}) \right \|_2 + L_{\rm{g}} \tilde{\beta}_{k} \|\tilde{\bm{v}}_{k}\|_2\\
        \leq & \left(\frac{\sigma_{\max}}{2} +\frac{L_{\rm{H}}}{2}   \right)\|\tilde{\bm{s}}_{k}\|_2^2 + \left\| \nabla m_{k}(\tilde{\bm{s}}_{k}) \right\|_2 \\
        & \quad + L_{\rm{g}} \alpha_2 \|\tilde{\bm{s}}_{k}\|_2^2 \cdot \|\tilde{\bm{v}}_{k}\|_2\\
        \leq & c_{5} \|\tilde{\bm{s}}_{k}\|_2^2 +  \left\| \nabla m_{k}(\tilde{\bm{s}}_{k}) \right\|_2,
    \end{split}
\end{equation*}
and
\begin{equation*}
    \begin{split}
        & \lambda_{\min}\left(\nabla^2 f(\bm{x}_{k+1})\right)\\
        \geq & \lambda_{\min}\left(\nabla^2 f(\bm{y}_{k+1})\right) - \left\| \nabla^2 f(\bm{x}_{k+1}) - \nabla^2 f(\bm{y}_{k+1}) \right\|_2\\
        \geq & \lambda_{\min}\left(\nabla^2 f(\bm{y}_{k+1})\right) - L_{\rm{H}} \left\| \bm{x}_{k+1} - \bm{y}_{k+1} \right\|_2\\
        \geq & \lambda_{\min}\left(\nabla^2 f(\bm{y}_{k+1})\right) - L_{\rm{H}} \tilde{\beta}_{k} \|\tilde{\bm{v}}_{k}\|_2 \\
        \geq & -\left( \frac{1}{2}\sigma_{\max} + L_{\rm{H}} \right) \left\| \tilde{\bm{s}}_{k} \right\|_2  - \frac{\sigma_{\max}}{2} \left|\|\bm{s}_{k}\|_2 - \|\tilde{\bm{s}}_{k}\|_2\right|\\
        & - L_{\rm{H}} \alpha_1 \|\tilde{\bm{s}}_{k}\|_2 \cdot \|\tilde{\bm{v}}_{k}\|_2\\
        \geq & - c_{6} \left\| \tilde{\bm{s}}_{k} \right\|_2  - \frac{\sigma_{\max}}{2} \left|\|\bm{s}_{k}\|_2 - \|\tilde{\bm{s}}_{k}\|_2\right|,
    \end{split}
\end{equation*}
where $c_{5}=\frac{1}{2}\max\left\{L_{\rm{H}} \gamma_1, \sigma_{\min}\right\} + \frac{1}{2} L_{\rm{H}} + \frac{\alpha_2 L_{\rm{g}}}{1-\tau} \left( \frac{24(f(\bm{x}_0 ) - f^{*})}{\eta_1 \sigma_{\min}} \right)^{1/3}$ and $c_{6}=\frac{1}{2}\sigma_{\max} + L_{\rm{H}} + \frac{\alpha_1 L_{\rm{H}}}{1-\tau} \left( \frac{24(f(\bm{x}_0 ) - f^{*})}{\eta_1 \sigma_{\min}} \right)^{1/3}$.

\section{Additional Experimental Results}
The trajectories of the empirical loss during the iterations are displayed in Figure \ref{fig_loss}, where we have found in Figure \ref{fig_grad} that the momentum in opposite direction helps the convergence of ARCm at the beginning of the algorithm.

\end{document}